\def\version{0.98}
\def\journal{LMP}
\def\titlep{Non-cocommutative C$^{*}$-bialgebra defined as the direct sum 
of free group C$^{*}$-algebras
\footnote{This is the revision of the previous version.}}
\documentclass[11pt]{article}
\usepackage{graphicx,ifthen,epic,eepic}
\usepackage{amssymb}
\usepackage{latexsym}

 at12pt

\newcommand{\qed}{\hbox{\rule[-2pt]{3pt}{6pt}}}
\newcommand{\qedh}{\hfill\qed \\}

\newcommand{\vv}{\vspace{.3in}}

\newcommand{\vep}{\varepsilon}

\newtheorem{Thm}{Theorem}[section]

\newtheorem{rem}[Thm]{Remark}

\newtheorem{defi}[Thm]{Definition}
\newtheorem{lem}[Thm]{Lemma}

\newtheorem{prop}[Thm]{Proposition}

\newtheorem{fact}[Thm]{Fact}
\newtheorem{fig}[Thm]{Figure}

\newcommand{\ww}{\vv\noindent}

\newcommand{\kn}{\Large\bf
$K\hspace{-.4cm} N$
\Large\bf\vv }

\def\cal#1{\mathcal #1}

\def\pr{{\it Proof.}\quad}

\def\co#1{{\cal O}_{#1}}
\def\brl{branching law}
\def\bfsnl{{\rm BFS}_{N}(\Lambda)}

\addtocounter{footnote}{1}
\def\sftt#1{
\setcounter{equation}{0}
\addtocounter{footnote}{1}
\section{#1}
}

\def\ssft#1{\subsection{#1}}
\def\sssft#1{\subsubsection{#1}}

\begin{document}
%
% Personal data
%
\def\autherp{Katsunori Kawamura}
\def\emailp{e-mail: kawamura@kurims.kyoto-u.ac.jp.}
\def\addressp{{\small {\it College of Science and Engineering, 
Ritsumeikan University,}}\\
{\small {\it 1-1-1 Noji Higashi, Kusatsu, Shiga 525-8577, Japan}}
}

\def\infw{\Lambda^{\frac{\infty}{2}}V}
\def\zhalfs{{\bf Z}+\frac{1}{2}}
\def\ems{\emptyset}
\def\pmvac{|{\rm vac}\!\!>\!\! _{\pm}}
\def\vac{|{\rm vac}\rangle _{+}}
\def\dvac{|{\rm vac}\rangle _{-}}
\def\ovac{|0\rangle}
\def\tovac{|\tilde{0}\rangle}
\def\expt#1{\langle #1\rangle}
\def\zph{{\bf Z}_{+/2}}
\def\zmh{{\bf Z}_{-/2}}
\def\brl{branching law}
\def\bfsnl{{\rm BFS}_{N}(\Lambda)}
\def\scm#1{S({\bf C}^{N})^{\otimes #1}}
\def\mqb{\{(M_{i},q_{i},B_{i})\}_{i=1}^{N}}
\def\zhalf{\mbox{${\bf Z}+\frac{1}{2}$}}
\def\zmha{\mbox{${\bf Z}_{\leq 0}-\frac{1}{2}$}}
\newcommand{\mline}{\noindent
\thicklines
\setlength{\unitlength}{.1mm}
\begin{picture}(1000,5)
\put(0,0){\line(1,0){1250}}
\end{picture}
\par
 }
\def\ptimes{\otimes_{\varphi}}
\def\delp{\Delta_{\varphi}}
\def\delps{\Delta_{\varphi^{*}}}
\def\gamp{\Gamma_{\varphi}}
\def\gamps{\Gamma_{\varphi^{*}}}
\def\sem{{\sf M}}
\def\sen{{\sf N}}
\def\hdelp{\hat{\Delta}_{\varphi}}
\def\tilco#1{\tilde{\co{#1}}}
\def\ndm#1{{\bf M}_{#1}(\{0,1\})}
\def\fs{{\cal F}{\cal S}({\bf N})}
\def\bigptimes{\bigotimes}
\def\bbn{{\Bbb N}}
\def\bbz{{\Bbb Z}}
\def\bbc{{\Bbb C}}
\def\bbr{{\Bbb R}}
\def\bbf{{\Bbb F}}
\def\bbt{{\Bbb T}}
\def\bbb{{\Bbb B}}
% Boldfont
\def\ba{\mbox{\boldmath$a$}}
\def\bb{\mbox{\boldmath$b$}}
\def\bc{\mbox{\boldmath$c$}}
\def\be{\mbox{\boldmath$e$}}
\def\bp{\mbox{\boldmath$p$}}
\def\bq{\mbox{\boldmath$q$}}
\def\bu{\mbox{\boldmath$u$}}
\def\bv{\mbox{\boldmath$v$}}
\def\bw{\mbox{\boldmath$w$}}
\def\bx{\mbox{\boldmath$x$}}
\def\by{\mbox{\boldmath$y$}}
\def\bz{\mbox{\boldmath$z$}}
\def\cdm#1{{\cal M}_{#1}(\{0,1\})}
%%%%%%%%%%%%%%%%%%%%%%%%%%%%%%%%%%%%%%%%%%%%%%%%
\def\titlepage{%\vspace{-4cm}

\noindent
{\bf 
\noindent
\thicklines
\setlength{\unitlength}{.1mm}
\begin{picture}(1000,0)(0,-300)
\put(0,0){\kn \knn\, for \journal\, Ver.\,\version}
\put(0,-50){\today,\quad {\rm file:} {\rm {\small \textsf{tit01.txt,\, J1.tex}}}}
\end{picture}
}
\vspace{-2.5cm}
\quad\\
{\small 
\footnote{
\begin{minipage}[t]{6in}
directory: \textsf{\fileplace}, \\
file: \textsf{\incfile},\, from \startdate
\end{minipage}
}}
\quad\\
\framebox{
%\noindent
\begin{tabular}{ll}
\textsf{Title:} &
\begin{minipage}[t]{4in}
\titlep
\end{minipage}
\\
\textsf{Author:} &\autherp
\end{tabular}
}
%\mline
{\footnotesize	
\tableofcontents }
}

%%%%%%%%%%%%%%%%%%%%%%%%%%%%%%%%%%%%%%%%%%%%%%%%

%%%%%%%%%%%%%%%%%%%%%%%%%%%%%%%%%%%%%%%%%%%%%%%%
\def\csrfn{C^{*}_{r}(\bbf_{n})}
\def\csrf#1{C^{*}_{r}(\bbf_{#1})}
\def\csfn{C^{*}(\bbf_{n})}
\def\csf#1{C^{*}(\bbf_{#1})}
\def\ul#1{\underline{#1}}
%
%%%%%%%%% Cut from here %%%%%%%%%%
%\input comm.txt
%%%%%%%%% End of Cut %%%%%%%%%
%
%
\setcounter{section}{0}
\setcounter{footnote}{0}
\setcounter{page}{1}
\pagestyle{plain}

%
% Title
%
%%%%%%%%%%%%%%%%%%%%%%%%%%%%%%%%%%%%%%%%%%
\title{\titlep}
\author{\autherp\thanks{\emailp}
\\
\addressp}
\date{}
\maketitle

%%%%%%%%%%%%%%%%%%%%%%%%%
%
% Abstract
%
\begin{abstract}
Let ${\Bbb F}_{n}$ be the free group of rank $n$ and 
let $\bigoplus C^{*}({\Bbb F}_{n})$ denote the direct sum of 
full group C$^{*}$-algebras $C^{*}({\Bbb F}_{n})$ of ${\Bbb F}_{n}$
$(1\leq n<\infty$).
We introduce a new comultiplication $\Delta_{\varphi}$ 
on $\bigoplus C^{*}({\Bbb F}_{n})$ 
such that $(\bigoplus C^{*}({\Bbb F}_{n}),\,\Delta_{\varphi})$
is a non-cocommutative C$^{*}$-bialgebra.
With respect to $\Delta_{\varphi}$,
the tensor product $\pi\otimes_{\varphi}\pi'$ of any two representations
$\pi$ and $\pi'$
of free groups is defined.
The operation $\ptimes$ is associative and non-commutative.
We compute its tensor product formulas 
of several representations.
\end{abstract}


%%%%%%%%%%%%%%%%%%%%%%%%%%%%%%
%
% Section 1
%
\sftt{Introduction}
\label{section:first}
A C$^{*}$-bialgebra is a generalization of bialgebra 
in the theory of C$^{*}$-algebras,
which was introduced in C$^{*}$-algebraic framework 
for quantum groups \cite{KV,MNW}. 
For example, if $G$ is a locally compact group,
then the full group C$^{*}$-algebra $C^{*}(G)$ of $G$ 
is a cocommutative C$^{*}$-bialgebra
with  respect to the standard (diagonal) comultiplication.

In this paper,
a C$^{*}$-bialgebra arising from certain group homomorphisms 
among free groups is given as follows:
Let ${\Bbb F}_{n}$ denote the free group of rank $n$ 
with free generators $g_{1}^{(n)},\ldots,g^{(n)}_{n}$.
For $n,m\geq 1$, 
define the group homomorphism $\phi_{n,m}$ from $\bbf_{nm}$
to $\bbf_{n}\times \bbf_{m}$ by
%
% Equation 1.1
%
\begin{equation}
\label{eqn:phib}
\phi_{n,m}(g_{m(i-1)+j}^{(nm)})
:= (g_{i}^{(n)},g_{j}^{(m)})
\quad(i=1,\ldots,n,\,j=1,\ldots,m).
\end{equation}
The map $\phi_{n,m}$ is well-defined 
on the whole of ${\Bbb F}_{nm}$
by the universality of $\bbf_{nm}$. 
Then the following diagram is commutative for each $n,m,l\geq 1$:

\noindent
%%%%%%%%%%%%%%%%%%%%%%%%%%%%%%%%%%%%%
%
\def\free{
\put(30,250){$\bbf_{nml}$}
\put(150,280){\vector(3,1){200}}
\put(150,350){$\phi_{n,ml}$}
\put(150,240){\vector(3,-1){200}}
\put(150,160){$\phi_{nm,l}$}
\put(420,350){$\bbf_{n}\times  \bbf_{ml}$}
\put(420,150){$\bbf_{nm}\times  \bbf_{l}$}
\put(650,350){\vector(3,-1){200}}
\put(700,350){$id_{n}\times  \phi_{m,l}$}
\put(650,170){\vector(3,1){200}}
\put(720,160){$\phi_{n,m}\times  id_{l}$}
\put(900,250){$\bbf_n\times  \bbf_m\times  \bbf_l$.}
}
\thicklines
%
% Figure 1.1
%
\begin{fig}
\label{fig:one}
\quad\\
\setlength{\unitlength}{.1mm}
\begin{picture}(1000,300)(-30,120)
\put(0,0){\free}
\end{picture}
\end{fig}

\noindent
Group homomorphisms in (\ref{eqn:phib})
can be lifted as 
$*$-homomorphisms $\varphi_{n,m}$ 
among full group C$^{*}$-algebras and their minimal tensors. 
For $\{\varphi_{n,m}\}$,
the following diagram is also commutative for each $n,m,l\geq 1$:

\noindent
%%%%%%%%%%%%%%%%%%%%%%%%%%%%%%%%%%%%%
%
\def\freeb{
\put(-50,250){$\csf{nml}$}
\put(150,280){\vector(2,1){150}}
\put(150,350){$\varphi_{n,ml}$}
\put(150,240){\vector(2,-1){150}}
\put(150,160){$\varphi_{nm,l}$}
\put(320,400){$\csf{n}\otimes \csf{ml}$}
\put(320,100){$\csf{nm}\otimes \csf{l}$}
\put(680,360){\vector(2,-1){150}}
\put(750,350){$id_{n}\otimes \varphi_{m,l}$}
\put(680,160){\vector(2,1){150}}
\put(750,160){$\varphi_{n,m}\otimes id_{l}$}
\put(850,250){$\csf{n}\otimes \csf{m}\otimes \csf{l}$.}
}
\thicklines
%
% Figure 1.2
%
\begin{fig}
\label{fig:two}
\quad\\
\setlength{\unitlength}{.1mm}
\begin{picture}(1000,400)(-30,60)
\put(0,0){\freeb}
\end{picture}
\end{fig}

By using $\{\varphi_{n,m}\}$,
we can construct a new comultiplication $\delp$
on the direct sum 
%
% Equation 1.2
%
\begin{equation}
\label{eqn:firstone}
\bigoplus C^{*}(\bbf_{n})
=C^{*}(\bbf_{1})\oplus C^{*}(\bbf_{2})
\oplus C^{*}(\bbf_{3})\oplus \cdots
\end{equation}
for all finite-rank free groups $\{\bbf_{n}:1\leq n<\infty\}$  
such that 
$(\bigoplus C^{*}(\bbf_{n}),\, \delp)$ 
is a non-cocommutative  C$^{*}$-bialgebra without antipode
(Theorem \ref{Thm:mainone}).

For any two unitary representations of free groups,
we can define the tensor product $\ptimes$ by using the
comultiplication $\delp$
which is {\it not} commutative
(Fact \ref{fact:noncommutative}).
Especially, the $\ptimes$-tensor product of any 
two quasi-regular representations is a direct sum of 
quasi-regular representations (Theorem \ref{Thm:decop}):
%
% Equation 1.3
%
\begin{equation}
\label{eqn:decopb}
\lambda_{{\Bbb F}_{n}/H'}\ptimes 
\lambda_{{\Bbb F}_{m}/H''}\cong 
\bigoplus_{\mu}
\lambda_{{\Bbb F}_{nm}/H_{\mu}}.
\end{equation}

In this section, we show our motivation, definitions and the main theorem.

%%%%%%%%%%%%%%%%%%%%%%%%%%%%%
%
% subsection 1.1
%  
\ssft{Motivation}
\label{subsection:firstone}
According to \cite{Herschend},
given two representations of a group $G$, 
their tensor product is a new representation
of $G$, which decomposes into a direct sum of indecomposable representations.
The problem of finding this decomposition is called 
the {\it Clebsch-Gordan problem}
and the resulting formula for the decomposition is called 
the {\it tensor product formula}
(or {\it Clebsch-Gordan formula} \cite{Herschend}).
A generalization of the Clebsch-Gordan problem for groups is to consider
modules over associative algebras instead of group algebras. However, there lies an
obvious obstruction in that there is no known way to define the tensor product of
two left modules over an arbitrary associative algebra. 
For group algebras, 
the extra structure coming from the group yields the tensor product.
For a bialgebra $A$, 
the associative tensor product of representations (=special modules) of $A$ can be
defined by using the comultiplicatoin. 
Hence one of most important motivations of 
the study of bialgebras is the tensor product of their representations.

We have studied a new kind of C$^{*}$-bialgebras which are defined as direct sums
of well-known C$^{*}$-algebras, for example,
Cuntz algebras,  UHF algebras, matrix algebras \cite{TS02} 
and Cuntz-Krieger algebras \cite{TS05}.
They are non-commutative and non-cocommutative,
and there never exist antipodes on them.
Such bialgebra structures do not appear before one takes direct sums. 
With respect to their comultiplications,
new tensor products among representations of these C$^{*}$-algebras 
and their tensor product formulas
were obtained
\cite{TS01,TS07}.
In \cite{TS02}, we gave a general method to construct  a C$^{*}$-bialgebra
from a given system of C$^{*}$-algebras and special $*$-homomorphisms among them.
The essential part of this construction
is how to construct such $*$-homomorphisms for each concrete example.
One of our interests is to construct new examples of C$^{*}$-bialgebra from
various C$^{*}$-algebras.

On the other hand, group C$^{*}$-algebras are important examples of C$^{*}$-algebras
\cite{BO,Davidson,Ped}.
Furthermore, quantum groups in the C$^{*}$-algebra approach are 
founded on the study of group C$^{*}$-algebras \cite{KV,MNW}.
%Of course, known examples of such quantum groups are 
%deformations of cocommutative C$^{*}$-bialgebras.
%(This is a well-known historical overview, but our study is not related to this history.)

Hence we consider to construct a new C$^{*}$-bialgebra associated 
with group C$^{*}$-algebras
by using a new comultiplication instead of their standard comultiplications.
In this paper, we choose free group C$^{*}$-algebras 
for this purpose,
and try to construct a new comultiplication on them
according to our method \cite{TS02}. 

%%%%%%%%%%%%%%%%%%%%%%%%%%
%
% subsection 1.2
%
\ssft{C$^{*}$-bialgebra}
\label{subsection:firsttwo}
In this subsection,
we review terminology about C$^{*}$-bialgebra according to \cite{ES,KV,MNW}.
For two C$^{*}$-algebras $A$ and $B$,
we write ${\rm Hom}(A,B)$ as the set of all $*$-homomorphisms 
from $A$ to $B$.
We assume that every tensor product $\otimes$ 
as below means the minimal C$^{*}$-tensor product.
%
% Definition 1.4
%
\begin{defi}
\label{defi:cstar}
A pair $(A,\Delta)$ is a C$^{*}$-bialgebra
if $A$ is a C$^{*}$-algebra and 
$\Delta\in {\rm Hom}(A,M(A\otimes A))$, 
where $M(A\otimes A)$ 
denotes the multiplier algebra of $A\otimes A$,
such that the linear span of $\{\Delta(a)(b\otimes c):a,b,c\in A\}$ 
is norm dense in $A\otimes A$ and 
the following holds:
%
% Equation 1.3
%
\begin{equation}
\label{eqn:bialgebratwo}
(\Delta\otimes id)\circ \Delta=(id\otimes\Delta)\circ \Delta.
\end{equation}
We call $\Delta$ the comultiplication of $A$.
\end{defi}

\noindent
We say that a C$^{*}$-bialgebra $(A,\Delta)$ is {\it strictly proper} 
if $\Delta(a)\in A\otimes A$ for any $a\in A$;
$(A,\Delta)$ is {\it unital}
if $A$ is unital and $\Delta$ is unital;
$(A,\Delta)$ is {\it counital}
if there exists $\vep\in {\rm Hom}(A,\bbc)$ such that
%
% Equation 1.4
%
\begin{equation}
\label{eqn:counit}
(\vep\otimes id)\circ \Delta= id = (id\otimes \vep)\circ\Delta.
\end{equation}
We call $\vep$ the {\it counit} of $A$ and write $(A,\Delta,\vep)$ 
as the counital C$^{*}$-bialgebra $(A,\Delta)$ with the counit $\vep$.
Remark that Definition \ref{defi:cstar} does not mean 
$\Delta(A)\subset A\otimes A$.
If $A$ is unital, then  $(A,\Delta)$ is strictly proper.
A {\it bialgebra} in the purely algebraic theory \cite{Abe,Kassel} means 
a unital counital strictly proper bialgebra with the unital counit
with respect to the algebraic tensor product,
which does not need to have an involution.
Hence a C$^{*}$-bialgebra is not a bialgebra in general.
In Definition \ref{defi:cstar},
if $A$ is unital and $\Delta $ is unital,
then the condition 
of the dense subspace in $A\otimes A$ can be omitted.

According to \cite{TS02},
we recall several notions of C$^{*}$-bialgebra.
%
% Definition 1.5
%
\begin{defi}
\label{defi:comodules}
\begin{enumerate}
%(i)
\item
For two C$^{*}$-bialgebras $(A_{1},\Delta_{1})$ and $(A_{2},\Delta_{2})$,
$f$ is a C$^{*}$-bialgebra morphism from 
$(A_{1},\Delta_{1})$ to $(A_{2},\Delta_{2})$ 
if $f$ is a non-degenerate $*$-homomorphism from $A_{1}$ to $M(A_{2})$
such that $(f\otimes f)\circ \Delta_{1}=\Delta_{2}\circ f$.
In addition, if $f(A_{1})\subset A_{2}$, then $f$ is called strictly proper.
%(ii)
\item
A map $f$ is a C$^{*}$-bialgebra endomorphism of
a C$^{*}$-bialgebra $(A,\Delta)$ if 
$f$ is a C$^{*}$-bialgebra morphism from $A$ to $A$.
In addition, if $f(A)\subset A$ and $f$ is bijective,
then $f$ is called a C$^{*}$-bialgebra automorphism of $(A,\Delta)$.
%(iii)
\item
A pair $(B,\Gamma)$ is a right comodule-C$^{*}$-algebra
of a C$^{*}$-bialgebra $(A,\Delta)$
if $B$ is a C$^{*}$-algebra and $\Gamma$ is a non-degenerate 
$*$-homomorphism from $B$ to $M(B\otimes A)$ such that
the following holds:
%
% Equation 1.5
%
\begin{equation}
\label{eqn:onethree}
(\Gamma\otimes id)\circ \Gamma=
(id\otimes \Delta)\circ \Gamma
\end{equation}
where both $\Gamma\otimes id$ and 
$id\otimes \Delta$ are extended to unital $*$-homomorphisms
from $M(B\otimes A)$ to $M(B\otimes A\otimes A)$.
The map $\Gamma$ is called the right coaction of $A$ on $B$.
%(iv)
\item
A proper C$^{*}$-bialgebra $(A,\Delta)$ satisfies the cancellation law 
if $\Delta(A)(I\otimes A)$ and $\Delta(A)(A\otimes I)$ are dense in $A\otimes A$
where $\Delta(A)(I\otimes A)$ and $\Delta(A)(A\otimes I)$ denote 
the linear spans of sets 
$\{\Delta(a)(I\otimes b):a,b\in A\}$ and $\{\Delta(a)(b\otimes I):a,b\in A\}$,
respectively. 
\end{enumerate}
\end{defi}

Let $(B,m,\eta,\Delta,\vep)$ be a bialgebra in the purely algebraic theory,
where $m$ is a multiplication and $\eta$ is a unit of the algebra $B$.
An endomorphism $S$ of $B$ is called an {\it antipode}
for $(B,m,\eta,\Delta,\vep)$ if $S$ satisfies
$m\circ (id\otimes S)\circ \Delta=\eta\circ \vep=
m\circ (S\otimes id)\circ \Delta$ \cite{Abe,Kassel}.

%%%%%%%%%%%%%%%%%%%%%%
%
% subsection 1.3
%
\ssft{Free group algebras and homomorphisms among them}
\label{subsection:firstthree}
In this subsection,
we briefly review free group C$^{*}$-algebras
\cite{BO,Davidson},
and introduce new homomorphisms among them
in order to define a comultiplication.

For $n=\infty, 1,2,3,\ldots$,
let $\bbf_{n}$ denote the free group of rank $n$
where we use the symbol ``$\infty$"  as the countable infinity for convenience
in this paper.
Let $({\cal K}_{n},\eta_{n})$ denote a direct sum of all irreducible representations 
(up to unitary equivalence) of 
the Banach algebra $\ell^{1}(\bbf_{n})$.
Let $\csfn$ denote the
{\it full group C$^{*}$-algebra} of $\bbf_{n}$,
which is defined as the C$^{*}$-algebra generated 
by the image of $\ell^{1}(\bbf_{n})$ by $\eta_{n}$.
Remark that $\csf{1}$ is $*$-isomorphic to 
the C$^{*}$-algebra
$C(\bbt)$ of all complex-valued continuous functions on
the torus $\bbt$.
With respect to the natural identification of 
the group algebra 
$\bbc\bbf_{n}$ over the coefficient field $\bbc$ with a subalgebra of $\csfn$,
$\bbc\bbf_{n}$ is dense in $\csfn$.
%There exists the tracial state on $\csfn$
%which is induced by the trivial representation of $\bbf_{n}$.
For $n=\infty,1,2,3,\ldots$,
let  $\{g_{i}^{(n)}\}$ be
the free generators of $\bbf_{n}$.
We also identify  $g_{i}^{(n)}$  with
the unitary $\eta_{n}(g_{i}^{(n)})$ in $\csfn$.

We introduce $*$-homomorphisms among $\csfn$'s as follows.
%
% Lemma 1.6
%
\begin{lem}
\label{lem:extension}
\begin{enumerate}
%(i)
\item
For $1\leq n,m<\infty$,
define the map $\varphi_{n,m}$ from $\csf{nm}$
to the minimal tensor product $\csf{n}\otimes  \csf{m}$ by
%
% Equation 1.7
%
\begin{equation}
\label{eqn:sseven}
\varphi_{n,m}(g_{m(i-1)+j}^{(nm)})
:=  g_{i}^{(n)}\otimes g_{j}^{(m)}
\quad(i=1,\ldots,n,\,j=1,\ldots,m).
\end{equation}
Then it is well-defined on the whole of  $\csf{nm}$ 
as a unital $*$-homomorphism.
%(ii)
\item
For $1\leq n<\infty$,
define
the map 
$\varphi_{\infty,n}$ 
from $\csf{\infty}$ to the minimal tensor product $\csf{\infty}\otimes \csfn$ by
%
% Equation 1.8
%
\begin{equation}
\label{eqn:sixteen}
\varphi_{\infty,n}(g^{(\infty)}_{n(i-1)+j})
:=  g^{(\infty)}_{i}\otimes g^{(n)}_{j}\quad(i\geq 1,\, j=1,\ldots,n).
\end{equation}
Then it is well-defined on the whole of  $\csf{\infty}$ 
as a unital $*$-homomorphism.
% (iii)
\item
If $n,m\geq 2$, then
$\varphi_{n,m}$ is not injective.
%(iv)
\item
Let $\csrf{n}$ denote the reduced group C$^{*}$-algebra of $\bbf_n$,
which is defined as the C$^{*}$-algebra 
generated by the image of the left regular representation of $\bbf_{n}$.
Then the map $\varphi_{n,m}$ in (\ref{eqn:sseven})
can not be extended as a $*$-homomorphism
from $\csrf{nm}$
to $\csrf{n}\otimes \csrf{m}$.
\end{enumerate}
\end{lem}

\noindent
Especially, $\varphi_{1,1}$
equals the standard comultiplication of $C^{*}(\bbf_{1})$.
The proof of Lemma \ref{lem:extension} will
 be given in $\S$ \ref{subsection:secondtwo}.

%%%%%%%%%%%%%%%%%%%%%%%%
%
% subsection 1.4
%
\ssft{Main theorem}
\label{subsection:firstfour}
In this subsection, we show our main theorem.
Let $\csfn$,  $\{g^{(n)}_{i}\}_{i=1}^{n}$,
$\bbc\bbf_{n}$,
 $\{\varphi_{n,m} \}_{n,m\geq 1}$
and 
 $\{\varphi_{\infty,n} \}_{n\geq 1}$
 be as in $\S$ \ref{subsection:firstthree}.
%
% Theorem 1.7
%
\begin{Thm}
\label{Thm:mainone}
Define the C$^{*}$-algebra ${\cal A}$ as the direct sum
%
% Equation 1.9
%
\begin{equation}
\label{eqn:onesix}
{\cal A}:=  \bigoplus_{1\leq n<\infty}\csfn
\end{equation}
and define 
$\delp\in {\rm Hom}({\cal A},{\cal A}\otimes {\cal A})$ 
and $\vep\in {\rm Hom}({\cal A},\bbc)$  by
%
% Equation 1.10
%
\begin{equation}
\label{eqn:oneseven}
\delp(x):=  \sum_{m,l;\, ml=n}\varphi_{m,l}(x)\quad
\mbox{when }x\in \csfn,
\end{equation}
%
% Equation 1.11
%
\begin{equation}
\label{eqn:coco}
\vep:=\vep_1\circ E_1  
\end{equation}
where 
$\varepsilon_1\in {\rm Hom}(C^*({\Bbb F}_1), {\Bbb C})$
%$\vep_{1}:\csf{1}\to \bbc$ 
is 
defined as $\vep_{1}|_{\bbf_{1}}=1$,
and $E_1$ is the projection from 
${\cal A}$ onto $\csf{1}$.
Then the following holds:
\begin{enumerate}
%(i)
\item
The C$^{*}$-algebra ${\cal A}$ is a strictly proper
counital C$^{*}$-bialgebra
with the comultiplication $\delp$ and the counit $\vep$.
%(ii)
\item
The C$^{*}$-bialgebra $({\cal A},\delp)$ satisfies the cancellation law.
%(iii)
\item
%The bialgebra $({\cal A},\delp,\vep)$ has no unit.
%But,
By the smallest unitization, $({\cal A},\delp,\vep)$ can be extended
to the unital counital C$^{*}$-bialgebra $(\tilde{{\cal A}},\hdelp,\tilde{\vep})$.
%(\cite{TS02}, Lemma 2.2).
%(iv)
\item
There never exists any antipode 
for any dense unital counital subbialgebra of $(\tilde{{\cal A}},\hdelp,\tilde{\vep})$ in (iii).
% (v)
\item
Define the algebraic direct sum 
$\bbc\bbf_{*}:=  
\oplus_{alg}\{\bbc\bbf_{n}: 1\leq n<\infty\}$.
Then 
$\delp(\bbc\bbf_{*})\subset 
\bbc\bbf_{*}\odot \bbc\bbf_{*}$
where $\odot$ means the algebraic tensor product,
and $\bbc\bbf_{*}$ is identified with
a $*$-subalgebra of ${\cal A}$
with respect to
the canonical embedding.
%(v)
\item
Define $\gamp \in
{\rm Hom}(\csf{\infty}, M(\csf{\infty}\otimes {\cal A}))$ by 
%
% Equation 1.12
%
\begin{equation}
\label{eqn:fifteen}
\gamp(x):=  \prod_{1\leq n<\infty}\varphi_{\infty,n}(x)
\quad(x\in \csf{\infty}),
\end{equation}
where
 we identify
the multiplier $M(\csf{\infty}\otimes {\cal A})$ with
the direct product $\prod_{n\geq 1}\csf{\infty}\otimes \csfn$.
Then  $\csf{\infty}$ is a 
right comodule-C$^{*}$-algebra of 
$({\cal A},\delp)$ with respect to the coaction $\gamp$.
\end{enumerate}
\end{Thm}

%
% Remark 1.8
% 
\begin{rem}
\label{rem:first}
{\rm
\begin{enumerate}
%(i)
\item
The R.H.S. in (\ref{eqn:oneseven})
is always a finite sum when $x\in \csfn$.
%(ii)
\item
The C$^{*}$-bialgebra
$({\cal A},\delp)$ is non-cocommutative.
In fact, the following holds:
%shows $\delp\ne \delp^{op}$:
%
% Equation 1.13
%
\begin{equation}
\label{eqn:utwo}
\delp(g_{2}^{(6)})
=
  g_{1}^{(1)}\otimes g_{2}^{(6)}
+g_{1}^{(2)}\otimes g_{2}^{(3)}
+g_{1}^{(3)}\otimes g_{2}^{(2)}
+g_{2}^{(6)}\otimes g_{1}^{(1)}.
\end{equation}
%(iii)
\item
In (\ref{eqn:onesix}),
every free group C$^{*}$-algebras $\csfn$ ($1\leq n<\infty$)
appear at once.
This is an essentially new structure of the class of free group C$^{*}$-algebras.
On the other hand,
$\csf{\infty}$ appears as a comodule-C$^{*}$-algebra
of $({\cal A},\,\delp)$.
This shows a certain naturality of this bialgebra structure.
\end{enumerate}
}
\end{rem}

In $\S$ \ref{section:second},
we prove Theorem \ref{Thm:mainone}.
In $\S$ \ref{section:third},
we show tensor product formulas 
of representations of $\bbf_{n}$'s
with respect to $\delp$, 
and show some C$^{*}$-bialgebra automorphisms.

%%%%%%%%%%%%%%%%%%%%%%%%%%%%%%%%%%%%%%%%%%%%%%%%%%%%%%%%
%
% Section 2
%
\sftt{Proofs of theorems}
\label{section:second}
In this  section, 
we prove Lemma \ref{lem:extension} and 
Theorem \ref{Thm:mainone}.
%%%%%%%%%%%%%%%%%%%%%%%%%%%%%%%%%%%%%%%
%
% Subsection 2.1
%
\ssft{C$^{*}$-weakly coassociative system}
\label{subsection:secondone}
According to $\S$ 3 in \cite{TS02},
we recall a general method to construct a C$^{*}$-bialgebra
from a set of C$^{*}$-algebras and $*$-homomorphisms among them.
A {\it monoid} is a set $\sem$ equipped with a binary associative operation 
$\sem\times \sem\ni(a,b)\mapsto ab\in \sem$,
and a unit with respect to the operation.
For example,
$\bbn=\{1,2,3,\ldots\}$ is an abelian monoid
with respect to the multiplication.
In order to show Theorem \ref{Thm:mainone},
we give a new definition of 
C$^{*}$-weakly coassociative system 
which is a generalization of the older in Definition 3.1 of \cite{TS02}.
%
% Definition 2.1
% 
\begin{defi}
\label{defi:axiom}
Let $\sem$ be a monoid with the unit $e$.
A data $\{(A_{a},\varphi_{a,b}):a,b\in \sem\}$
is a C$^{*}$-weakly coassociative system (= C$^{*}$-WCS) over $\sem$ if 
$A_{a}$ is a unital C$^{*}$-algebra for $a\in \sem$
and $\varphi_{a,b}$ is a unital $*$-homomorphism
from $A_{ab}$ to $A_{a}\otimes A_{b}$
for $a,b\in \sem$ such that
\begin{enumerate}
%(i)
\item
for all $a,b,c\in \sem$, the following holds:
%
% Equation 2.1
%
\begin{equation}
\label{eqn:wcs}
(id_{a}\otimes \varphi_{b,c})\circ \varphi_{a,bc}
=(\varphi_{a,b}\otimes id_{c})\circ \varphi_{ab,c}
\end{equation}
where $id_{x}$ denotes the identity map on $A_{x}$ for $x=a,c$,
%(ii)
\item
there exists a counit $\vep_{e}$ of $A_{e}$ 
such that $(A_{e},\varphi_{e,e},\vep_{e})$ 
is a counital C$^{*}$-bialgebra,
%(iii)
\item
for each $a\in\sem$, the following holds:
%
% Equation 2.2
%
\begin{equation}
\label{eqn:new}
(\vep_{e}\otimes id_{a})\circ \varphi_{e,a}=id_{a}
=(id_{a}\otimes \vep_{e})\circ \varphi_{a,e}.
\end{equation}
\end{enumerate}
\end{defi}

\noindent 
The condition (\ref{eqn:new}) is weaker than the older,
``$\varphi_{e,a}(x)=I_{e}\otimes x$ and
$\varphi_{a,e}(x)=x\otimes I_{e}$ for $x\in A_{a}$ and $a\in \sem$"
(\cite{TS02}, Definition 3.1) .
In fact, the older definition satisfies (\ref{eqn:new}).
From the new definition, the same result holds as follows.
%
% Theorem 2.2
% 
\begin{Thm}
\label{Thm:mainthree}
(\cite{TS02}, Theorem 3.1).
Let $\{(A_{a},\varphi_{a,b}):a,b\in \sem\}$ be a C$^{*}$-WCS 
over a monoid $\sem$.
Assume that $\sem$ satisfies that 
%
% Equation 2.3
%
\begin{equation}
\label{eqn:finiteness}
\#{\cal N}_{a}<\infty \mbox{ for each }a\in \sem
\end{equation}
where ${\cal N}_{a}:= \{(b,c)\in \sem\times \sem:\,bc=a\}$.
Define C$^{*}$-algebras 
\[A_{*}:=   \oplus \{A_{a}:a\in \sem\},\quad
C_{a}:=  
\oplus \{A_{b}\otimes A_{c}:(b,c)\in {\cal N}_{a}\}
\quad (a\in\sem). \]
Define $\Delta^{(a)}_{\varphi}\in{\rm Hom}(A_{a},C_{a})$,
$\Delta_{\varphi}
\in {\rm Hom}(A_{*}, A_{*}\otimes A_{*})$ and
$\vep\in {\rm Hom}(A_{*},\bbc)$ by 
\[\Delta^{(a)}_{\varphi}(x):=  \sum_{(b,c)\in {\cal N}_{a}}
\varphi_{b,c}(x)\quad(x\in A_{a}),\quad 
\Delta_{\varphi}:=  \oplus\{\Delta_{\varphi}^{(a)}:a\in \sem\},\]
%
% Equation 2.4
%
\begin{equation}
\label{eqn:counittwo}
\vep:=\vep_e\circ E_e  
\end{equation}
where $E_e$ denotes the projection
from $A_{*}$ onto $A_e$.
Then $(A_{*},\delp,\vep)$ 
is a strictly proper counital C$^{*}$-bialgebra.
\end{Thm}
%
% Proof
%
\pr
By (\ref{eqn:finiteness}), 
$\delp^{(a)}$ is well-defined.
Furthermore,
$C_{a}$ is unital 
and $\Delta^{(a)}_{\varphi}$ is unital for each $a$.
Since $\sem\times \sem=\coprod_{a\in \sem}{\cal N}_{a}$,
$A_{*}\otimes A_{*}=\oplus\{A_{f}\otimes A_{g}:f,g\in \sem\}
=\oplus \{C_{a}:a\in \sem\}$.
Since $\delp^{(a)}$ is unital for each $a$,
$\delp$ is non-degenerate.
From (\ref{eqn:wcs}), the following holds for $x\in A_{a}$:
%
% Equation 2.5
%
\begin{equation}
\label{eqn:appfour}
\begin{array}{rl}
\{(\Delta_{\varphi}\otimes id)\circ \Delta_{\varphi}\}(x)
=&
\sum_{b,c,d\in \sem,\,bcd=a}(\varphi_{b,c}\otimes id_{d})(\varphi_{bc,d}(x))\\
=&
\sum_{b,c,d\in \sem,\,bcd=a}(id_{b}\otimes \varphi_{c,d})(\varphi_{b,cd}(x))\\
=&\{(id\otimes \Delta_{\varphi})\circ \Delta_{\varphi}\}(x).
\end{array}
\end{equation}
Hence
$(\Delta_{\varphi}\otimes id)\circ \Delta_{\varphi}
=(id\otimes \Delta_{\varphi})\circ \Delta_{\varphi}$ on $A_{*}$.
Therefore $\Delta_{\varphi}$ is a comultiplication of $A_{*}$.
On the other hand,
for $x\in A_{a}$, we see that
%
% Equation 2.6
%
\begin{equation}
\label{eqn:arrayone}
\begin{array}{rl}
\{(\vep\otimes id)\circ \delp\}(x)=&(\vep\otimes id)(\delp^{(a)}(x))\\
=&\sum_{(b,c)\in {\cal N}_{a}}(\vep\otimes id)(\varphi_{b,c}(x))\\
=&(\vep_{e}\otimes id_{a})(\varphi_{e,a}(x))\\
=& x\qquad(\mbox{from }(\ref{eqn:new})).
\end{array}
\end{equation}
Hence $(\vep\otimes id)\circ \delp= id$.
In like wise, we see that $(id\otimes \vep)\circ \delp= id$.
Therefore $\vep$ is a counit of $(A_{*},\delp)$.
In consequence, we see that $(A_{*},\Delta_{\varphi},\vep)$ 
is a counital C$^{*}$-bialgebra.
By definition, $(A_{*},\Delta_{\varphi})$ is strictly proper.
\qedh

\noindent
We call $(A_{*},\Delta_{\varphi},\vep)$ in 
Theorem \ref{Thm:mainthree} by a (counital)
{\it C$^{*}$-bialgebra} associated with 
$\{(A_{a},\varphi_{a,b}):a,b\in \sem\}$.

The following lemma holds independently of the generalization 
in Definition \ref{defi:axiom}(iii).
%
% Lemma 2.3
%
\begin{lem}
\label{lem:unitization}
For the following C$^{*}$-WCS $\{(A_{a},\varphi_{a,b}):a,b\in \sem\}$, 
we assume the condition (\ref{eqn:finiteness}).
\begin{enumerate}
%(i)
\item
(\cite{TS02},  Lemma 2.2)
For a given strictly proper non-unital counital C$^{*}$-bialgebra
$(A,\Delta,\vep)$,
let $\tilde{A}:=  A\oplus \bbc$ denote the
smallest unitization of $A$.
Then there exist a unique extension
$(\hat{\Delta},\tilde{\vep})$
of $(\Delta,\vep)$ on $\tilde{A}$
such that  $(\tilde{A},\hat{\Delta},\tilde{\vep})$
is a strictly proper unital counital C$^{*}$-bialgebra.
%(ii)
\item
(\cite{TS02}, Lemma 3.2)
For a C$^{*}$-WCS $\{(A_{a},\varphi_{a,b}):a,b\in \sem\}$
over $\sem$,
let $(A_{*},\delp,\vep)$ be as in Theorem \ref{Thm:mainthree}
and let $(\tilde{A}_{*},\hdelp,\tilde{\vep})$ be the smallest unitization
of $(A_{*},\delp,\vep)$ in (i). 
Assume that any element in $\sem$ 
has no left inverse except the unit $e$.
Then the antipode for any dense unital counital subbialgebra of 
$(\tilde{A}_{*},\hdelp,\tilde{\vep})$ never exists.
%(iii)
\item
(\cite{TS02}, Lemma 3.1)
Let $\{(A_{a},\varphi_{a,b}):a,b\in \sem\}$ be a 
C$^{*}$-WCS over a monoid $\sem$
and let $(A_{*},\delp)$ be as in Theorem \ref{Thm:mainthree}
associated with $\{(A_{a},\varphi_{a,b}):a,b\in \sem\}$. 
Define
%
% Equation 2.7
%
\begin{equation}
\label{eqn:notation}
X_{a,b}:=  \varphi_{a,b}(A_{ab})(A_{a}\otimes I_{b}),\quad
Y_{a,b}:=  \varphi_{a,b}(A_{ab})(I_{a}\otimes A_{b})\quad(a,b\in \sem)
\end{equation}
where
$\varphi_{a,b}(A_{ab})(A_{a}\otimes I_{b})$
and $\varphi_{a,b}(A_{ab})(I_{a}\otimes A_{b})$
mean
the linear spans 
of 
$\{\varphi_{a,b}(x)(y\otimes I_{b}):x\in A_{ab},\,y\in A_{a}\}$
and
$\{\varphi_{a,b}(x)(I_{a}\otimes y):x\in A_{ab},\,y\in A_{b}\}$,
respectively.
If both $X_{a,b}$ and $Y_{a,b}$ 
are dense in $A_{a}\otimes A_{b}$ for each $a,b\in \sem$,
then $(A_{*},\delp)$ satisfies the cancellation law.
%(iv)
\item
(\cite{TS02}, Theorem 3.2)
For a C$^{*}$-WCS $\{(A_{a},\varphi_{a,b}):a,b\in \sem\}$ 
over a monoid $\sem$,
assume that $B$ is a unital C$^{*}$-algebra and 
a set $\{\varphi_{B,a}:a\in \sem\}$ of unital $*$-homomorphisms
such that $\varphi_{B,a}\in {\rm Hom}(B,B\otimes A_{a})$ for each $a\in \sem$ 
and the following holds:
%
% Equation 2.8
%
\begin{equation}
\label{eqn:comodulefour}
(\varphi_{B,a}\otimes id_{b})\circ \varphi_{B,b}
=(id_{B}\otimes \varphi_{a,b})\circ \varphi_{B,ab}\quad (a,b\in \sem).
\end{equation}
Then $B$ is a right comodule-C$^{*}$-algebra of 
the C$^{*}$-bialgebra $(A_{*},\delp)$ with the unital coaction
$\gamp:=  \prod_{a\in\sem}\varphi_{B,a}$.
\end{enumerate}
\end{lem}

%%%%%%%%%%%%%%%%%%%%%%%%%%%%%%%
%
% subsection 2.2
%
\ssft{Homomorphisms among free groups}
\label{subsection:secondtwo}
In this subsection, we show
properties of  $\phi_{n,m}$ in (\ref{eqn:phib})
and prove Lemma \ref{lem:extension}.
%
% Lemma 2.4
%
\begin{lem}
\label{lem:representation}
For each $n\geq 1$,
we write $1$ as the unit of $\bbf_{n}$.
\begin{enumerate}
%(i)
\item
For any 
$x\in \bbf_{n}$,
there exists
$(y,z)\in \bbf_{m}\times \bbf_{nm}$
such that $\phi_{n,m}(z)=(x,y)$.
%(ii)
\item
For any 
$y\in \bbf_{m}$,
there exists
$(x,z)\in \bbf_{n}\times \bbf_{nm}$
such that $\phi_{n,m}(z)=(x,y)$.
%(iii)
\item
For any $(x,y)\in \bbf_{n}\times \bbf_{m}$,
there exists $(x',z)\in \bbf_{n}\times \bbf_{nm}$
such that $\phi_{n,m}(z)(x',1)=(x,y)$.
%(iv)
\item
For any $(x,y)\in \bbf_{n}\times \bbf_{m}$,
there exists $(y',z)\in \bbf_{m}\times \bbf_{nm}$
such that $\phi_{n,m}(z)(1,y')=(x,y)$.
%(v)
\item
When $n,m\geq 2$,
$\phi_{n,m}$ is not injective. 
\end{enumerate}
\end{lem}
%
% Proof
%
\pr
(i)
Let
 $a_{1},\ldots,a_{n}$,
 $b_{1},\ldots,b_{m}$,
 $c_{1},\ldots,c_{nm}$ be the free generators of 
$\bbf_{n},
\bbf_{m},
\bbf_{nm}$, respectively.
Assume that $x\in\bbf_{n}$ is written as a reduced word 
$x=a_{i_{1}}^{\vep_1}\cdots  a_{i_{l}}^{\vep_l}$
where $\vep_{i}=1$ or $-1$ for $i=1,\ldots,l$.
For  example, define $(y,z)\in \bbf_{m}\times \bbf_{nm}$ by
$y:=  b_{1}^{\vep_{1}}\cdots b_{1}^{\vep_{l}}$
and 
$z:=  
c_{m(i_{1}-1)+1}^{\vep_{1}}\cdots 
c_{m(i_{l}-1)+1}^{\vep_{l}}$.
Then $y$ belongs to the abelian subgroup generated by
the single element $b_{1}$,
and it
is not always a reduced word in $\bbf_{m}$.
Then the statement holds for $(y,z)$.

\noindent
(ii) As the proof of (i), this is proved.

\noindent
(iii)
From (ii),
we can find $(x'',z)\in \bbf_{n}\times \bbf_{nm}$
such that 
$\phi_{n,m}(z)=(x'',y)$.
Define $x':=  (x'')^{-1}x$,
then the statement holds.

\noindent
(iv) As the proof of (iii), this is proved from (i).

\noindent
(v)
Let $c_{1},\ldots,c_{nm}$ be as in the proof of (i).
For $i,l\in \{1,\ldots,n\},\,
k,j\in\{1,\ldots,m\}$,
define $x(i,l;j,k)\in\bbf_{nm}$ by
%
% Equation 2.9
%
\begin{equation}
\label{eqn:ften}
x(i,l;j,k):=  
c_{m(i-1)+j}\,
c_{m(i-1)+k}^{-1}\,
c_{m(l-1)+k}\,
c_{m(l-1)+j}^{-1}.
\end{equation}
Then $x(i,l;j,k)\ne 1$ when $k\ne j,\, i\ne l$,
but $x(i,l;j,k)\in \ker \phi_{n,m}$ for any
$i,l,j,k$.
\qedh

\noindent
In the proof of Lemma \ref{lem:representation}(v),
if $n=m=2$,  then
the reduced word $c_{1}c_{2}^{-1}c_{4}c_{3}^{-1}$ in $\bbf_{4}$ satisfies
$\phi_{2,2}(c_{1}c_{2}^{-1}c_{4}c_{3}^{-1})=(1,1)$.\\

\noindent
{\it Proof of Lemma \ref{lem:extension}}
(i)
Let 
$\phi_{n,m}$ be as in (\ref{eqn:phib}) and let
$({\cal K}_{n},\eta_{n})$ be as in $\S$ \ref{subsection:firstthree}.
Define the unitary representation $\varphi_{n,m}^{0}$
of $\bbf_{nm}$ on ${\cal K}_{n}\otimes {\cal K}_{m}$
by $\varphi_{n,m}^{0}:=  
(\eta_{n}\otimes \eta_{m})\circ \phi_{n,m}$.
The representation
$\varphi_{n,m}^{0}$ is well-defined by the universality of $\bbf_{nm}$.
Since the image of $\varphi_{n,m}^{0}$ is included 
in $\csfn\otimes \csf{m}$,
$\varphi_{n,m}^{0}$ is uniquely extended
to  $\varphi_{n,m}$ in (\ref{eqn:sseven})
such that 
$\varphi_{n,m}(\eta_{nm}(x))=\varphi_{n,m}^{0}(x)$
for each $x\in \bbf_{nm}$
(\cite{BO}, Proposition 2.5.2).
Hence the statement holds.

\noindent
(ii)
In analogy with (i), the statement holds.

\noindent
(iii)
For $x(i,l;j,k)$ in (\ref{eqn:ften}),
we see that $\varphi_{n,m}(x(i,l;j,k)-1)=0$ for each $i,l,j,k$.
Hence the statement holds.

\noindent
(iv)
If such an extension $\tilde{\varphi}_{n,m}$  of $\varphi_{n,m}$
exists,
then $\tilde{\varphi}_{n,m}$ must be injective because $\csrf{nm}$ is simple
when $nm\geq 2$
(\cite{Davidson}, Corollary VII.7.5 and its proof).
On the other hand,
$\tilde{\varphi}_{n,m}$ never be injective  from (iii).
\qedh

%%%%%%%%%%%%%%%%%%%%%%%%%%%%
%
% subsection 2.3
%
\ssft{Proof of Theorem \ref{Thm:mainone}}
\label{subsection:secondthree}
We prove  Theorem \ref{Thm:mainone} in this subsection.
Let $\bbn:=\{1,2,3,\ldots\}$.
Remark that (\ref{eqn:finiteness}) holds
for any element in the multiplicative monoid $(\bbn,\cdot)$.

\noindent
(i)
From Theorem \ref{Thm:mainthree},
it is sufficient to show that
$\{(\csfn,\varphi_{n,m}):n,m\in \bbn\}$
is a C$^{*}$-WCS over the monoid $\bbn$.
By the definition of $\varphi_{n,m}$ in (\ref{eqn:sseven}),
we can verify that
$(\varphi_{n,m}\otimes id_{l})
\circ \varphi_{nm,l}
=
(id_{n}\otimes \varphi_{m,l})
\circ \varphi_{n,ml}$ for $n,m,l\in \bbn$
where $id_{a}$ denotes the identity map on $\csf{a}$ for $a=n,l$.
Hence (\ref{eqn:wcs}) is satisfied.
On the other hand,
since $\vep_{1}|_{\bbf_{1}}=1$,
$\{(\vep_{1}\otimes id_{n})
\circ \varphi_{1,n}\}(g^{(n)}_{j})
=(\vep_{1}\otimes id_{n})(g^{(1)}_{1}\otimes g^{(n)}_{j})
=\vep_{1}(g^{(1)}_{1})\,g^{(n)}_{j}=g^{(n)}_{j}$
for each $j=1,\ldots,n$ and $n\in \bbn$.
By the same token,
we obtain $(id_{n}\otimes \vep_{1})\circ \varphi_{n,1}=id_{n}$.
Hence (\ref{eqn:new}) is verified.
Therefore $\{(\csfn,\varphi_{n,m}):n,m\in \bbn\}$
is a C$^{*}$-WCS over the monoid $\bbn$.

\noindent
(ii)
%{\it Proof of cancellation law.}
For $n,m\in \bbn$,
define  three subsets 
${\cal P}_{n,m},{\cal Q}_{n,m},{\cal R}_{n,m}$ 
of $\csfn\otimes \csf{m}$
by
%
% Equation 2.10-12
%
\begin{eqnarray}
\label{eqn:triplea}
{\cal P}_{n,m}:=  &
\{\varphi_{n,m}(z)(x\otimes I_{m})
:x\in\bbf_{n},\,
z\in\bbf_{nm}
\},
%\varphi_{n,m}(\bbf_{nm}) (\bbf_{n}\otimes I_{m}),
\\ \nonumber
\\
\label{eqn:tripleb}
{\cal Q}_{n,m}
:=  &
\{\varphi_{n,m}(z)(I_{n}\otimes y)
:y\in\bbf_{m},\,
z\in\bbf_{nm}\},
%\varphi_{n,m}(\bbf_{nm})
%(I_{n}\otimes \bbf_{m}),
\\ \nonumber
\\
\label{eqn:triplec}
{\cal R}_{n,m}:=  &
\{x\otimes y
:
x\in\bbf_{n},\,
y\in\bbf_{m}
\}.
%\bbf_{n}\otimes \bbf_{m}.
\end{eqnarray}
Then their linear spans are dense subspaces
of $\varphi_{n,m}(\csf{nm})(\csfn\otimes I_{m})$, 
$\varphi_{n,m}(\csf{nm})(I_{n}\otimes \csf{m})$ and
$\csfn\otimes \csf{m}$,  respectively.
From Lemma \ref{lem:representation}(iii),
it is sufficient to show that 
${\cal R}_{n,m}\subset {\cal P}_{n,m}$
and ${\cal R}_{n,m}\subset {\cal Q}_{n,m}$.

We prove ${\cal R}_{n,m}\subset {\cal P}_{n,m}$ as follows:
For $(x,y)\in \bbf_{n}\times \bbf_{m}$,
there exists $(x',z)\in \bbf_{n}\times \bbf_{nm}$
such that $\phi_{n,m}(z)(x',1)=(x,y)$
from Lemma \ref{lem:representation}(iii).
By definitions of $\phi_{n,m}$ and $\varphi_{n,m}$,
this implies 
$\varphi_{n,m}(z)(x'\otimes I_{m})=x\otimes y$.
Therefore ${\cal R}_{n,m}\subset {\cal P}_{n,m}$.

In a similar fashion,
we obtain ${\cal R}_{n,m}\subset {\cal Q}_{n,m}$
from Lemma \ref{lem:representation}(iv).
Hence the statement holds.

\noindent
(iii)
From (i) and Lemma \ref{lem:unitization}(i), the statement holds.

\noindent
(iv)
Remark that the monoid $\bbn$ 
has no left invertible element except the unit $1$.
From the proof of (i) and Lemma \ref{lem:unitization}(ii), the statement holds.

\noindent
(v)
From the proof of (i) and
the definition of $\varphi_{n,m}$ in (\ref{eqn:sseven}),
we see that
$\varphi_{n,m}(\bbc\bbf_{nm})\subset \bbc\bbf_{n}\odot \bbc\bbf_{m}$.
This implies the statement.

\noindent
(vi)
By definition, we see that 
$(\varphi_{\infty,n}\otimes id_{m})\circ \varphi_{\infty,m}
=(id_{\infty}\otimes \varphi_{n,m})\circ \varphi_{\infty,nm}$
for $n,m\in\bbn$.
From Lemma \ref{lem:unitization}(iv)
for $\varphi_{\csf{\infty},n}:=  \varphi_{\infty,n}$,
the statement holds.
\qedh

%%%%%%%%%%%%%%%%%%%%%%%%%%%
%
% Section 3
%
\sftt{Tensor product formulas of representations, and automorphisms}
\label{section:third}
In this section,
we show tensor product formulas
of unitary representations of $\bbf_{n}$'s
with respect to the comultiplication $\delp$ in Theorem \ref{Thm:mainone},
and C$^{*}$-bialgebra automorphisms.

%%%%%%%%%%%%%%%%%%%%%%%%%%%%%%%%%%%%%
%
%  subsection 3.1
% 
\ssft{General facts about representations}
\label{subsection:thirdone}
We introduce a new tensor product 
of representations of ${\Bbb F}_n$'s
and show tensor product formulas of 
quasi-regular representations.
%%%%%%%%%%%%%%%%%%%%%%%%%%%%%%%%%%%%%%%%%%%%
%
% subsubsection 3.1.1
%
\sssft{Representations of $\bbf_{n}$}
\label{subsubsection:thirdoneone}
We identify $\bbf_{n}$ with the unitary subgroup of 
$\csfn$ with respect to the canonical embedding.
Let ${\rm Rep}_{u}\bbf_{n}$ denote the class 
of all unitary representations of $\bbf_{n}$.
For $(\pi,\pi')\in {\rm Rep}_{u}\bbf_{n}\times {\rm Rep}_{u}\bbf_{m}$,
define the new representation $\pi\ptimes \pi'\in {\rm Rep}_{u}\bbf_{nm}$ by
%
% Equation 3.1
%
\begin{equation}
\label{eqn:tone}
\pi\ptimes \pi':=  
(\pi\otimes \pi')\circ \varphi_{n,m},
\end{equation}
where $\varphi_{n,m}$ is as in (\ref{eqn:sseven}).
Then we see that the new operation $\ptimes $ is associative, and 
it is distributive with respect to the direct sum.
Furthermore, $\ptimes$ is well-defined 
on the unitary equivalence classes of representations.
It will be shown that $\ptimes$ is non-commutative
in Fact \ref{fact:noncommutative}.

%%%%%%%%%%%%%%%%%%%
% 
% subsubsection  3.1.2
% 
\sssft{Tensor product formulas of 1-dimensional unitary
representations of free groups}
\label{subsubsection:thirdonetwo}
In this subsection,
we show tensor product formulas of 1-dimensional unitary
representations of free groups with respect to
$\ptimes$ in (\ref{eqn:tone})
as a basic example of tensor product formula.
For groups $G$ and $H$,
let ${\rm Hom}(G,H)$ denote
the set of all homomorphisms from $G$ to $H$.
Define 
%
% Equation 3.2
%
\begin{equation}
\label{eqn:ch}
{\rm Ch}({\Bbb F}_n):={\rm Hom}({\Bbb F}_n,U(1)).
\end{equation}
Since an element in ${\rm Ch}({\Bbb F}_n)$ can be regarded 
as a one-dimensional unitary representation of ${\Bbb F}_n$,
we can regard $\pi_1\ptimes \pi_2$ as an element 
in ${\rm Ch}({\Bbb F}_{nm})$ for 
$\pi_1\in {\rm Ch}({\Bbb F}_n)$ and 
$\pi_2\in {\rm Ch}({\Bbb F}_m)$.

Let $g_1,\ldots,g_n$ be free generators of ${\Bbb F}_n$.
For $z=(z_1,\ldots,z_n)\in {\Bbb T}^n:=U(1)^{n}$,
define
$\chi_z^{(n)}\in {\rm Ch}({\Bbb F}_n)$ by
%
% Equation 3.3
%
\begin{equation}
\label{eqn:chiz}
\chi_z^{(n)}(g_i):=z_i\quad(i=1,\ldots,n).
\end{equation}
Clearly, $\chi_z^{(n)}$ is irreducible
for any $z\in {\Bbb T}^n$.
For $z,w\in {\Bbb T}^n$,
$\chi_z^{(n)}$ and $\chi_w^{(n)}$ are unitarily equivalent
if and only if $z=w$.
By the correspondence 
%
% Equation 3.4
%
\begin{equation}
\label{eqn:correspondence}
{\Bbb T}^n\ni z\mapsto \chi^{(n)}_z\in
{\rm Ch}({\Bbb F}_n),
\end{equation}
we see that ${\rm Ch}({\Bbb F}_n)$ is equivalent to ${\Bbb T}^{n}$
as a set.

For $z\in {\Bbb T}^n$
and $w\in{\Bbb T}^m$,
define the Kronecker product 
$z\boxtimes w\in {\Bbb T}^{nm}$
as 
$(z\boxtimes w)_{m(i-1)+j}:=z_iw_j$ for
$(i,j)\in\{1,\ldots,n\}\times\{1,\ldots,m\}$.
By definition, the following holds:
%
% Equation 3.5
%
\begin{equation}
\label{eqn:chizn}
\chi_z^{(n)}\ptimes \chi_w^{(m)}=\chi_{z\boxtimes w}^{(nm)}
\quad (z\in {\Bbb T}^n,\,w\in{\Bbb T}^m).
\end{equation}
This implies that 
the correspondence in (\ref{eqn:correspondence})
gives a semigroup isomorphism from 
$(\bigcup_{n\geq 1}{\rm Ch}({\Bbb F}_n),\,\ptimes)$
to $(\bigcup_{n\geq 1}{\Bbb T}^n,\,\boxtimes)$.
This shows a naturality of the operation $\ptimes$.

%
% Fact 3.1
%
\begin{fact}
\label{fact:noncommutative}
The operation $\ptimes$ in (\ref{eqn:tone})
is non-commutative as 
the following sense:
There exist $n,m\geq 2$ and representations $\pi_1,\pi_2$ 
of ${\Bbb F}_n$ and
${\Bbb F}_m$, respectively such that 
$\pi_1\ptimes \pi_2$ and 
$\pi_2\ptimes \pi_1$ are not unitarily equivalent.
\end{fact}
%
% Proof
%
\pr
In (\ref{eqn:chizn}),
let $(n,m)=(2,3)$
and $z=(1,-1)\in{\Bbb T}^2$ and $w=(1,1,1)\in{\Bbb T}^3$.
Then
$\chi_z^{(2)}\ptimes \chi_w^{(3)}
=\chi_{z\boxtimes w}^{(6)}
\not\cong 
\chi_{w\boxtimes z}^{(6)}
= \chi_w^{(3)}\ptimes\chi_z^{(2)}$
because $z\boxtimes w=(1,1,1,-1,-1,-1)$
and $w\boxtimes z=(1,-1,1,-1,1,-1)$.
\qedh

%%%%%%%%%%%%%%%%%%%%%%%%%%%%%%%%%%%%%%%%%%%%
%
% subsubsection 3.1.3
%
\sssft{Quasi-regular representations}
\label{subsubsection:thirdonethree}
In this subsection,
we review quasi-regular representations
of discrete groups,
and show the general formula of the $\ptimes$-tensor product of quasi-regular 
representations of free groups.

For a discrete group $\Gamma$
and a subgroup $\Gamma_0$,
let $\Gamma/\Gamma_0$ denote the left coset space,
that is, $\Gamma/\Gamma_0:=\{x\Gamma_0:x\in\Gamma\}$.
Define the {\it quasi-regular representation} \cite{Hatem}
(or {\it permutation representation} \cite{Mackey})
 $(\ell^{2}(\Gamma/\Gamma_0),\lambda_{\Gamma/\Gamma_0})$  
of $\Gamma$ associated with $\Gamma_0$
as the natural (unitary) left action
of $\Gamma$ on the standard basis of $\ell^{2}(\Gamma/\Gamma_0)$:
%
% Equation 3.6
%
\begin{equation}
\label{eqn:cosetb}
\lambda_{\Gamma/\Gamma_0}:\Gamma\curvearrowright  \ell^{2}(\Gamma/\Gamma_0).
\end{equation}
Especially,
the regular representation $\lambda$
and the trivial representation ${\bf 1}$
of $\Gamma$ are quasi-regular representations associated with
subgroups $\{e\}$ and $\Gamma$ of $\Gamma$, respectively.
For the trivial representation ${\bf 1}_{\Gamma_0}$ of $\Gamma_0$,
$\lambda_{\Gamma/\Gamma_0}$
coincides with the induced representation
${\rm Ind}_{\Gamma_0}^{\Gamma}({\bf 1}_{\Gamma_0})$.

Since $\Gamma/\Gamma_0$ is a $\Gamma$-homogeneous space,
$\lambda_{\Gamma/\Gamma_0}$ is a cyclic representation.
When $\Gamma$ acts on a set $X$,
the permutation representation of $\Gamma$
on $\ell^2(X)$ are decomposed into the direct sum of 
quasi-regular representations as follows:
%
% Equation 3.7
%
\begin{equation}
\label{eqn:orbits}
\ell^2(X)\cong \bigoplus_{\mu}\ell^2(\Gamma/H_{\mu})
\end{equation}
where
$H_{\mu}$ is a subgroup of $\Gamma$ such that
$X_{\mu}\cong \Gamma/H_{\mu}$ for 
the orbit decomposition $ X=\coprod_{\mu}X_{\mu}$
with respect to the $\Gamma$-action. 
About the irreducibility and unitary equivalence
of quasi-regular representations,
see Appendix \ref{section:appone}.

Next,
we consider the tensor product $\ptimes$ in (\ref{eqn:tone})
for quasi-regular representations of free groups.
%
% Theorem  3.2
%
\begin{Thm}
\label{Thm:decop}
Let $H',H''$ be  subgroups of ${\Bbb F}_n$
and ${\Bbb F}_m$, respectively.
For $\phi_{n,m}$ in (\ref{eqn:phib}),
define the left action $\tilde{\phi}_{n,m}$ of ${\Bbb F}_{nm}$
on the direct product set $X:={\Bbb F}_{n}/H'\times {\Bbb F}_{m}/H''$ by 
%
% Equation 3.8
%
\begin{equation}
\label{eqn:tilde}
\tilde{\phi}_{n,m}(g)(xH',yH''):=(g'xH',g''yH'')
\quad (\,(xH',yH'')\in X,\,g\in {\Bbb F}_{nm}\,)
\end{equation}
where $(g',g''):=\phi_{n,m}(g)$.
With respect to the action $\tilde{\phi}_{n,m}$,
let $X=\coprod_{\mu}X_{\mu}$ be the orbit decomposition
and choose $H_{\mu}$ as a stabilizer subgroup of ${\Bbb F}_{nm}$
associated with $X_{\mu}$.
%where $X_{\mu}:=\tilde{\phi}_{n,m}({\Bbb F}_{nm})x_{\mu}$
%for some $x_{\mu}\in X$,
Then the following holds:
%
% Equation 3.9
%
\begin{equation}
\label{eqn:decop}
\lambda_{{\Bbb F}_{n}/H'}\ptimes 
\lambda_{{\Bbb F}_{m}/H''}\,\cong \,
\bigoplus_{\mu}
\lambda_{{\Bbb F}_{nm}/H_{\mu}}.
\end{equation}
\end{Thm}
%
% Proof
%
\pr
Let $\{\xi'_a:a\in {\Bbb F}_{n}/H'\}$, $\{\xi''_b:b\in {\Bbb F}_{nm}/H''\}$
denote standard bases of 
$\ell^2({\Bbb F}_{n}/H')$
and 
$\ell^2({\Bbb F}_{nm}/H'')$,
respectively. 
By definition,
$(\lambda_{{\Bbb F}_{n}/H'}\ptimes 
\lambda_{{\Bbb F}_{m}/H''})(g)(\xi'_a\otimes \xi_b'')
=
\xi'_{g'a}\otimes \xi_{g''b}''$ for $g\in{\Bbb F}_{nm})$
where $(g',g''):=\phi_{n,m}(g)$.
On the other hand,
the permutation representation $(\ell^2(X),L)$ of ${\Bbb F}_{nm}$
by $\tilde{\phi}_{n,m}$ satisfies
$L_g\xi_{(a,b)}
=\xi_{\tilde{\phi}_{n,m}(g)(a,b)}
=\xi_{(g'a,\,g''b)}$
where $\{\xi_{(a,b)}:(a,b)\in X\}$ denotes the standard basis of $\ell^2(X)$.
By the natural identification
of $\ell^2(X)$ with $\ell^2({\Bbb F}_{n}/H')\otimes 
\ell^2({\Bbb F}_{m}/H'')$,
we see that $L$ and $\lambda_{{\Bbb F}_{n}/H'}\ptimes 
\lambda_{{\Bbb F}_{m}/H''}$
are unitarily equivalent.

By definition,
$X_{\mu}\cong {\Bbb F}_{nm}/H_{\mu}$ 
as ${\Bbb F}_{nm}$-homogeneous spaces,
and 
the statement holds from 
the decomposition in (\ref{eqn:orbits}).
\qedh

\noindent
Theorem \ref{Thm:decop} states
that the $\ptimes$-tensor product of any two quasi-regular representations
is decomposed into the direct sum of quasi-regular representations.
That is, the category of direct sums of quasi-regular representations
of free groups
is closed with respect to the $\ptimes$-tensor product.
This shows a naturality of the $\ptimes$-tensor product.
In $\S$ \ref{subsection:thirdtwo} and \ref{subsection:thirdthree},
we will show
concrete examples of 
the formula (\ref{eqn:decop}). 

%%%%%%%%%%%%%%%%%%%%%%%%%%%%%%%%%%%%%%%%%%%%%%
%
% subsection 3.2
%
\ssft{Tensor product of some irreducible representations}
\label{subsection:thirdtwo}
In this subsection, we show tensor product formulas of some irreducible
quasi-regular representations  of $\bbf_{n}$'s
as examples of Theorem \ref{Thm:decop}.

We review some irreducible quasi-regular representations of $\bbf_{n}$
\cite{AWW,Yoshizawa}.
Let $g_{1},\ldots,g_{n}$ be the free generators of $\bbf_{n}$.
Fix $i\in\{1,\ldots,n\}$ and
let $H_{i}^{(n)}$ denote the abelian subgroup of $\bbf_{n}$ generated
by the single element $g_{i}$:
%
% Equation 3.10
%
\begin{equation}
\label{eqn:hin}
H_{i}^{(n)}:=\{g_i^{l}:l\in{\Bbb Z}\}\subset {\Bbb F}_n.
\end{equation}
%
% Proposition 3.3
% 
\begin{prop}
\label{prop:irreducible}
\begin{enumerate}
%(i)
\item
For any $i=1,\ldots,n$,
$\lambda_{{\Bbb F}_n/H_{i}^{(n)}}$ is irreducible.
%(ii)
\item
$
\lambda_{{\Bbb F}_n/H_{i}^{(n)}}
\cong 
\lambda_{{\Bbb F}_n/H_{j}^{(n)}}$
if and only if $i=j$.
\end{enumerate}
\end{prop}
%
% Proof
%
\pr
See Appendix \ref{section:appone}.
\qedh

We show the tensor product formula of 
$\lambda_{{\Bbb F}_n/H^{(n)}_{i}}$'s
in Proposition \ref{prop:irreducible}.
For the map $\phi_{n,m}$ in (\ref{eqn:phib}),
define the subgroup $G_{n,m}$ of ${\Bbb F}_{nm}$ by
%
% Equation 3.11
%
\begin{equation}
\label{eqn:gnm}
G_{n,m}:=  \ker \phi_{n,m}.
\end{equation}
From Lemma \ref{lem:representation}(v),
$G_{n,m}\ne\{1\}$ when $n,m\geq 2$ and 
$G_{n,1}= G_{1,n}= \{1\}$ for any $n\geq 1$.
Since $G_{n,m}$ is a normal subgroup of ${\Bbb F}_{nm}$,
$G_{n,m}H:=\{gh:(g,h)\in G_{n,m}\times H\}$ is also a subgroup
of ${\Bbb F}_{nm}$ and $G_{n,m}H=HG_{n,m}$
for any subgroup $H$ of ${\Bbb F}_{nm}$.

%
% Theorem 3.4
%
\begin{Thm}
\label{Thm:ldotss}
Let $H^{(n)}_{i}$ and $G_{n,m}$ 
be as in (\ref{eqn:hin}) and (\ref{eqn:gnm}),
respectively.
Define $K_{n,m,l}:=G_{n,m}H^{(nm)}_{l}$
for $l=1,\ldots,nm$.
\begin{enumerate}
%(i)
\item
For $n,m\geq 1$
and $(i,j)\in\{1,\ldots,n\}\times \{1,\ldots,m\}$,
%
% Equation 3.12
%
\begin{equation}
\label{eqn:lambda}
\lambda_{{\Bbb F}_n/H_{i}^{(n)}}\ptimes 
\lambda_{{\Bbb F}_m/H_{j}^{(m)}}
\cong \lambda_{{\Bbb F}_{nm}/K_{n,m,m(i-1)+j}}.
\end{equation}
%(ii)
\item
For any $l=1,\ldots,nm$,
$\lambda_{{\Bbb F}_{nm}/K_{n,m,l}}$
is irreducible.
%(iii)
\item
$\lambda_{{\Bbb F}_{nm}/K_{n,m,l}}$
and 
$\lambda_{{\Bbb F}_{nm}/K_{n,m,l'}}$
are unitarily equivalent
if and only if $l=l'$.
\end{enumerate}
\end{Thm}
%
% Proofs 
% 
\pr
See Appendix \ref{section:apptwo}.
\qedh

%%%%%%%%%%%%%%%%%%%%%%%%%%%%%%%%%%%%%%%%%%%%
%
% subsection 3.3
%
\ssft{Tensor product of regular representations}
\label{subsection:thirdthree}
In this subsection,
we consider regular representations of $\bbf_{n}$'s
as examples of Theorem \ref{Thm:decop}.

We recall characterization of representations by using 
von Neumann algebras.
A nondegenerate representation $\pi$ of a C$^{*}$-algebra ${\cal A}$
is said to be {\it (pure) type {\rm X}}
if $\pi({\cal A})^{''}$ is of type X 
for X$=$I, II, III, II$_{1}$, II$_{\infty}$ \cite{Blackadar2006}.
% p315
For a group $G$ and unitary representation $U$ of $G$,
$U$ is said to be {\it type {\rm X}}
if $U(G)''$ is of type X 
for X$=$I, II, III, II$_{1}$, II$_{\infty}$.
$U$ is {\it factor} if the center of $U(G)''$ is trivial.

In the proof of Theorem \ref{Thm:decop},
$\lambda_{{\Bbb F}_n}=\lambda_{{\Bbb F}_{n}/H'}$
and
$\lambda_{{\Bbb F}_m}=\lambda_{{\Bbb F}_{m}/H''}$
when
$H'=\{1\}\subset {\Bbb F}_n$ and $H''=\{1\}\subset {\Bbb F}_{m}$.
%
% Proposition 3.5
%
\begin{prop}
\label{prop:lambdab}
Let $n,m\geq 1$.
\begin{enumerate}
%(i)
\item
$\lambda_{{\Bbb F}_n}\ptimes \lambda_{{\Bbb F}_m}\cong 
(\lambda_{{\Bbb F}_{nm}/G_{n,m}})^{\oplus \infty}$
where $G_{n,m}$ is as in (\ref{eqn:gnm}).
%(ii)
\item
If $n,m\geq 2$, then
$\lambda_{{\Bbb F}_{nm}/G_{n,m}}$
is a type {\rm II}$_1$ factor representation.
\end{enumerate}
\end{prop}
%
% Proof
%
\pr
(i)
Let $X:={\Bbb F}_n\times {\Bbb F}_m$
and let $\tilde{\phi}_{n,m}$ be the action of ${\Bbb F}_{nm}$
on $X$ in (\ref{eqn:tilde}).
From Lemma \ref{lem:representation}(iii), 
we see that the orbit decomposition is given as follows:
%
% Equation 3.13
%
\begin{equation}
\label{eqn:xcoprod}
X=\coprod_{x\in{\Bbb F}_n}
{\cal O}_x,\quad
{\cal O}_x:=\{
 \tilde{\phi}_{n,m}(y)(x,1):y\in {\Bbb F}_{nm}\}.
\end{equation}
Furthermore,
the equivalence
$({\cal O}_x,\tilde{\phi}_{n,m}|_{{\cal O}_x})
\cong ({\Bbb F}_{nm}/G_{n,m},L)$
holds 
as ${\Bbb F}_{nm}$-homogeneous spaces 
for all $x\in {\Bbb F}_n$
where $L$ denotes the natural left action of ${\Bbb F}_{nm}$
on ${\Bbb F}_{nm}/G_{n,m}$.
Therefore the equivalence 
$(X,\tilde{\phi}_{n,m})\cong 
({\Bbb F}_{nm}/G_{n,m},L)^{\#{\Bbb F}_{n}}$
of ${\Bbb F}_{nm}$-homogeneous spaces holds.
From this and the proof of Theorem \ref{Thm:decop},
the statement holds.

\noindent
(ii) See Appendix \ref{section:appthree}.
\qedh

%%%%%%%%%%%%%%%%%%%%%%%%%%%%%%%%%
%
% subsection 3.4
%
\ssft{Automorphisms}
\label{subsection:thirdfour}
In this subsection,
we show examples of some C$^{*}$-bialgebra automorphism
of $({\cal A},\delp)$.
For $t\in \bbr$,
define $\alpha_{t}^{(n)}\in {\rm Aut}\csfn$ 
by 
%
% Equation 3.14
%
\begin{equation}
\label{eqn:alpha}
\alpha_{t}^{(n)}(g_{i}^{(n)}):=   e^{\sqrt{-1}t\log n}\,g_{i}^{(n)}
\quad (i=1,\ldots,n).
\end{equation}
Then
$\alpha^{(*)}_{t}:=  \oplus_{n\geq 1}\alpha^{(n)}_{t}$
is a C$^{*}$-bialgebra automorphism
of $({\cal A},\delp)$
such that 
$\alpha^{(*)}_{t}\circ \alpha^{(*)}_{s}=\alpha^{(*)}_{t+s}$
for $s,t\in\bbr$.

Define $\beta^{(n)}\in {\rm Aut}\csfn$
by 
%
% Equation 3.15
%
\begin{equation}
\label{eqn:beta}
\beta^{(n)}(g_{i}^{(n)}):=  g_{n-i+1}^{(n)}
\quad (i=1,\ldots,n).
\end{equation}
Then
$\beta^{(*)}:=  \oplus_{n\geq 1}\beta^{(n)}$
is a C$^{*}$-bialgebra automorphism of
$({\cal A},\delp)$
such that $\beta^{(*)}\circ \beta^{(*)}=id$.

The automorphism $\beta^{(*)}$ commutes $\alpha_{t}^{(*)}$ 
for each $t$.
Hence these give the action of the group $\bbr\times (\bbz/2\bbz)$
on the C$^{*}$-bialgebra $({\cal A},\delp)$.

\ww
{\bf Acknowledgment:}
The author would like to express his sincere thanks to
reviewers for correcting errors in the previous version.

%%%%%%%%%%%%%%%%%%%%%%%%%%%%%%%%%%%%%%%%%%%%%%
\appendix

\section*{Appendix}

%%%%%%%%%%%%
%
% Appendix A
%
\sftt{Applications of commensurator to quasi-regular representations
of discrete groups}
\label{section:appone}
Recall that two subgroups $\Gamma_0$ and $\Gamma_1$ 
of a group $\Gamma$ are {\it commensurable} if 
$\Gamma_0\cap  \Gamma_1$
is of finite index in both $\Gamma_0$ and $\Gamma_1$,
that is, $[\Gamma_0,\Gamma_0\cap \Gamma_1]\cdot
[\Gamma_1,\Gamma_0\cap \Gamma_1]<\infty$,
in such case,
we write $\Gamma_0\approx \Gamma_1$.
Define the {\it commensurator} ${\rm Com}_{\Gamma}(\Gamma_0)$ 
of $\Gamma_0$ in $\Gamma$ as
%
% Equation A.1
%
\begin{equation}
\label{eqn:comg}
{\rm Com}_{\Gamma}(\Gamma_0):=\{g\in\Gamma:\Gamma_0\approx
g\Gamma_0 g^{-1}\}.
\end{equation}
Then the inclusions 
$\Gamma_0\subset {\rm Com}_{\Gamma}(\Gamma_0)\subset \Gamma$
of subgroups hold.
According to \cite{BH},
we review applications of commensurator to quasi-regular representations
of discrete groups by Mackey \cite{Mackey}.
% 3.f ?
%
% Theorem A.1
%
\begin{Thm}(\cite{BH}, Theorem 2.1)
\label{Thm:mackey}
Let $\Gamma$ be a discrete group and let $\Gamma_0,\Gamma_1$
be subgroups of $\Gamma$.
\begin{enumerate}
%(i)
\item
The representation $(\ell^2(\Gamma/\Gamma_0),\lambda_{\Gamma/\Gamma_0})$
of $\Gamma$
is irreducible if and only if ${\rm Com}_{\Gamma}(\Gamma_0)=\Gamma_0$.
%(ii)
\item
Assume ${\rm Com}_{\Gamma}(\Gamma_i)=\Gamma_i$ for $i=0,1$.
Then 
$\lambda_{\Gamma/\Gamma_0}$ and $\lambda_{\Gamma/\Gamma_1}$
are unitarily equivalent
if and only if  $\Gamma_0$ and $\Gamma_1$
are quasiconjugate in $\Gamma$, that is,
there exists $g\in\Gamma$ such that $\Gamma_0\approx g\Gamma_1 g^{-1}$.
\end{enumerate}
\end{Thm}

%
% Lemma A.2
%
\begin{lem}
\label{lem:multiple}
Let $H_{i}^{(n)}$ be as in (\ref{eqn:hin}).
\begin{enumerate}
%(i)
\item
If $a,b\in{\Bbb F}_n$ satisfy $ab=ba$, then there exists $w\in {\Bbb F}_n$
such that $a=w^l$ and $b=w^{l'}$
for some $l,l'\in {\Bbb Z}$.
%(ii)
\item
If $a,b\in{\Bbb F}_n$ satisfy $a^lb^k=b^ka^l$ 
for some $l,k\in {\Bbb Z}$,
then $ab=ba$.
%(ii)
\item
If $a\in{\Bbb F}_n$ satisfies $ab=ba$
for some $b\in H^{(n)}_i$ and $b\ne 1$,
then $a\in H^{(n)}_i$.
%(iii)
\item
If $i,j\in\{1,\ldots,n\}$ and 
$g\in{\Bbb F}_n$ satisfy
$H^{(n)}_i\approx gH_{j}^{(n)}g^{-1}$,
then $i=j$ and $g\in H_{i}^{(n)}$.
\end{enumerate}
\end{lem}
%
% Proof
%
\pr
Let $c_1,\ldots,c_n$ be
free generators of ${\Bbb F}_n$
and let $H_i:=H^{(n)}_i$.
\\
(i)
See \cite{MKS}, p42, 6.

\noindent
(ii)
See \cite{MKS}, p41, 4.

\noindent
(iii)
By (i), both $a$ and $b$ can be written as
$a=w^{l}$
and $b=w^{l'}$ for some $w\in{\Bbb F}_n$
and $l,l'\in {\Bbb Z}$.
By the choice of $b$,
$b=c_i^{k}$ for some $k\in {\Bbb Z}$
and $k\ne 0$. 
Hence $c_i^{k}=w^{l'}$.
Therefore $c_i^kw^{l'}=w^{l'}c_i^k$.
From (ii), $c_i w=wc_i$.
From (i), we see $w\in H_i$.
Therefore $a=w^{l}\in H_i$.

\noindent
(iv)
If $i,j$ and $g$ satisfy the assumption,
then $[H_i, \,
H_i\cap  gH_{j}g^{-1}]<\infty$ by definition.
Since $\#H_{i}=\infty$,
there exists $x\in 
H_i\cap  gH_{j}g^{-1}$ such that $x\ne 1$.
Then $c_i^l=x=gc_j^{l'}g^{-1}$
for some $l,l'\in {\Bbb Z}$.
For $\chi_z^{(n)}$ in (\ref{eqn:chiz}),
$z_i^{l}=\chi_z^{(n)}(c_i^l)
=\chi_z^{(n)}(gc_j^{l'}g^{-1})
=\chi_z^{(n)}(c_j^{l'})=z_j^{l'}$
for any $z\in {\Bbb T}^n$.
This implies $i=j$ and $l=l'$.
In consequence,
$c_i^l=gc_i^{l}g^{-1}$
for some $l\in {\Bbb Z}$.
By the choice of $x$ and $l$, $l\ne 0$.
From (iii), $g\in H_{i}$.
\qedh

\noindent
{\it Proof of Proposition \ref{prop:irreducible}.}
(i) 
From Lemma \ref{lem:multiple}(iii),
${\rm Com}_{{\Bbb F}_{n}}(H_i^{(n)})=H_i^{(n)}$.
By Theorem \ref{Thm:mackey}(i), the statement holds.

\noindent
(ii)
From Lemma \ref{lem:multiple}(iii),
$H_i^{(n)}$ and 
$H_j^{(n)}$ are quasiconjugate if and only if $i=j$.
From this and 
Theorem \ref{Thm:mackey}(ii),
the statement holds.
\qedh

%%%%%%%%%%%%%%%%%%%%%%%%%%%%%
%
% Section B
%
\sftt{Proof of Theorem \ref{Thm:ldotss}}
\label{section:apptwo}
In this section, we prove Theorem \ref{Thm:ldotss}.
Let $\phi_{n,m},H^{(n)}_{i},G_{n,m}$
be as in (\ref{eqn:phib}), (\ref{eqn:hin}) and (\ref{eqn:gnm}), respectively.
Since $G_{n,m}$ is a normal subgroup of ${\Bbb F}_{nm}$,
we can define the quotient group 
%
% Equation B.1
%
\begin{equation}
\label{eqn:qnm}
Q_{n,m}:={\Bbb F}_{nm}/G_{n,m}. 
\end{equation}
Then $Q_{1,n}\cong Q_{n,1}\cong {\Bbb F}_n$ for any $n\geq 1$.

Let $g_{1},\ldots,g_n$ be
free generators of ${\Bbb F}_n$.
Define $p^{(n)}\in{\rm Hom}({\Bbb F}_n,{\Bbb Z})$ by
%
% Equation B.2
%
\begin{equation}
\label{eqn:pg}
p^{(n)}(g):=
\vep_1+\cdots+\vep_k\quad  \mbox{when } 
g=g_{i_{1}}^{\vep_1}\cdots g_{i_{k}}^{\vep_k}\in{\Bbb F}_n\\
\end{equation}
where $\vep_i\in\{1,-1\}$ for $i=1,\ldots,k$.
By definition, the following holds.
%
% Fact B.1
%
\begin{fact}
\label{fact:pg}
\begin{enumerate}
%(i)
\item
For $g\in{\Bbb F}_{nm}$, 
if $\phi_{n,m}(g)=(g',g'')$, then $p^{(nm)}(g)=p^{(n)}(g')=p^{(m)}(g'')$.
%(ii)
\item
If $g\in G_{n,m}$,
then $p^{(nm)}(g)=0$.
%(iii)
\item
For any $i=1,\ldots,n$,
the restriction 
$p^{(n)}|_{H^{(n)}_{i}}:H^{(n)}_{i}\to {\Bbb Z}$
is an isomorphism.
%(iv)
\item
Let $\hat{p}^{(nm)}:Q_{n,m}\to {\Bbb Z}$ by
$\hat{p}^{(nm)}(\hat{g}):=p^{(nm)}(g)$ for $\hat{g}=gG_{n,m}\in Q_{n,m}$
and define the subgroup $\hat{H}^{(nm)}_l$
of $Q_{n,m}$ by
%
% Equation B.3
%
\begin{equation}
\label{eqn:hath}
\hat{H}^{(nm)}_{l}:=\{hG_{n,m}\in Q_{n,m}:h\in H_l^{(nm)}\}.
\end{equation}
Then $\hat{p}$ is well-defined and is a group homomorphism.
Furthermore,
the restriction $\hat{p}^{(nm)}|_{\hat{H}_l^{(nm)}}:\hat{H}_l^{(nm)}\to {\Bbb Z}$
is an isomorphism
for any $l=1,\ldots,nm$.
\end{enumerate}
\end{fact}

\noindent
Remark that $\hat{H}_l^{(1\cdot n)}\cong \hat{H}_l^{(n\cdot 1)}\cong H_l^{(n)}$
for any $n\geq 1$.

%
% Lemma B.2
%
\begin{lem}
\label{lem:hhnm}
Let $n,m\geq 1$,$(i,j)\in\{1,\ldots,n\}\times
\{1,\ldots,m\}$ and $g\in {\Bbb F}_{nm}$.
\begin{enumerate}
%(i)
\item
If $\phi_{n,m}(g)\in 
H^{(n)}_{i}\times
H^{(m)}_j$,
then there exists $h\in H^{(nm)}_{m(i-1)+j}$
such that $\phi_{n,m}(h)=\phi_{n,m}(g)$.
%(ii)
\item
If $\phi_{n,m}(g)\in 
H^{(n)}_{i}\times
H^{(m)}_j$,
then
$g\in G_{n,m}H^{(nm)}_{m(i-1)+j}$.
\end{enumerate}
\end{lem}
%
% Proof
%
\pr
(i)
Let $(g',g''):=\phi_{n,m}(g)$.
From Fact \ref{fact:pg}(i),
$p^{(n)}(g')=p^{(nm)}(g)=p^{(m)}(g'')$.
When $l:=p^{(nm)}(g)$,
$g'=(g_{i}^{(n)})^l$
and 
$g''=(g_{j}^{(m)})^l$ by Fact \ref{fact:pg}(iii).
Hence 
$h:=(g^{(nm)}_{m(i-1)+j})^l\in H^{(nm)}_{m(i-1)+j}$ satisfies the relation.

\noindent
(ii)
From (i),
$\phi_{n,m}(h^{-1}g)=(1,1)$ for some $h\in H^{(nm)}_{m(i-1)+j}$.
Therefore
$h^{-1}g\in G_{n,m}$
and $g\in hG_{n,m}\subset H^{(nm)}_{m(i-1)+j}G_{n,m}=G_{n,m}H^{(nm)}_{m(i-1)+j}$.
\qedh

%
% Lemma B.3
%
\begin{lem}
\label{lem:equiv}
For $n,m\geq 1$ and $(i,j)\in\{1,\ldots,n\}\times
\{1,\ldots,m\}$,
two ${\Bbb F}_{nm}$-homogeneous spaces
$({\Bbb F}_{nm}/(G_{n,m}H^{(nm)}_{m(i-1)+j}),\,L)$
and $({\Bbb F}_{n}/H_{i}^{(n)}\times {\Bbb F}_{m}/H_{j}^{(m)},
\,\tilde{\phi}_{n,m})$
are equivalent
where $L$ denotes the natural left action of ${\Bbb F}_{nm}$ on 
${\Bbb F}_{nm}/K$ and 
$\tilde{\phi}_{n,m}$ is as in (\ref{eqn:tilde}).
\end{lem}
%
% Proof
%
\pr
Rewrite 
$K:=G_{n,m}H^{(nm)}_{m(i-1)+j}$,
$\phi:=\phi_{n,m}$,
$H':=H^{(n)}_{i}$ and $H'':=H^{(m)}_{j}$ here.
Define the map
%
% Equation  B.4
%
\begin{equation}
\label{eqn:theta}
\theta:{\Bbb F}_{nm}/K
\to 
{\Bbb F}_n/H'\times 
{\Bbb F}_m/H'';\quad
\theta([x]):=([x'],[x''])
\end{equation}
where $(x',x''):=\phi(x)$
and $[a]$ denotes the coset
with the representative $a$. 
By definition, we see that
the map $\theta$ is well-defined and satisfies
$\theta \circ L_g=\tilde{\phi}_{n,m}(g)\circ \theta$
for all $g\in{\Bbb F}_{nm}$.
It is sufficient to show that 
$\theta$ is injective and surjective.

\noindent
(i) Injectivity:
For $[x],[y]\in {\Bbb F}_{nm}/K$,
assume  $\theta([x])=\theta([y])$.
Then
$[x']=[y']$
and $[x'']=[y'']$.
Hence
$x'=y'h'$
and $x''=y''h''$
for some $(h',h'')\in H'\times H''$.
Therefore $\phi(x)=(x',x'')=(y'h',y''h'')
=\phi(y)(h',h'')$.
Hence
$\phi(y^{-1}x)=(h',h'')$.
From Lemma \ref{lem:hhnm}(ii),
$y^{-1}x\in K$.
Therefore $[x]=[y]$.
Hence $\theta$ is injective.

\noindent
(ii) Surjectivity:
For $([x'],[x''])\in 
{\Bbb F}_n/H'\times 
{\Bbb F}_m/H''$,
there exists
$(w',z)\in {\Bbb F}_{n}\times {\Bbb F}_{nm}$
such that
$\phi(z)(w',1)=(x',x'')$
from Lemma \ref{lem:representation}(iii).
Let $g_1^{(n)},\ldots,g_n^{(n)}$ be 
free generators of ${\Bbb F}_n$.
Assume 
$w'=(g_{j_1}^{(n)})^{\vep_1}\cdots (g_{j_l}^{(n)})^{\vep_l}$
for $\vep_i\in\{1,-1\}$.
Let $h'':=(g^{(m)}_{j})^{\vep_1}\cdots (g^{(m)}_{j})^{\vep_l}
=(g^{(m)}_{j})^{\vep_1+\cdots +\vep_l}$.
Then
$(x',x''h'')=\phi(z)(w',h'')=\phi(z)\phi(w)=\phi(zw)$
for some $w\in {\Bbb F}_{nm}$.
Therefore
$([x'],[x''])=([x'],[x''h''])=\theta([zw])$.
Hence $\theta$ is surjective.
\qedh

Let $\chi_{z}^{(n)}$ be as in (\ref{eqn:chiz}) and 
let ${\Bbb T}^n\boxtimes {\Bbb T}^m:=\{z\boxtimes w:(z,w)\in 
{\Bbb T}^n\times {\Bbb T}^m\}$.
For $u\in {\Bbb T}^n\boxtimes {\Bbb T}^m$,
we see $G_{n,m}\subset \ker\chi_{u}^{(nm)}$
by (\ref{eqn:chizn}).
From this,
$\chi_{u}^{(nm)}(gh)=\chi_{u}^{(nm)}(g)$
for any $g\in {\Bbb F}_{nm}$
and $h\in G_{n,m}$.
Hence
%
% Equation B.5
%
\begin{equation}
\label{eqn:hatchi}
\hat{\chi}_{u}^{(nm)}:Q_{n,m}\to U(1);
\quad \hat{\chi}_{u}^{(nm)}(gG_{n,m}):=\chi_{u}^{(nm)}(g)
\end{equation}
is well-defined for any $u\in 
{\Bbb T}^n\boxtimes {\Bbb T}^m$ as a homomorphism.

%
% Lemma B.5
%
\begin{lem}
\label{lem:hij}
Let $n,m\geq 1$ and 
let $\hat{H}^{(nm)}_{l}$ be as in (\ref{eqn:hath}).
If $i,j\in\{1,\ldots,nm\}$ and 
$\hat{g}\in Q_{n,m}$ satisfy
$\hat{H}^{(nm)}_i\approx \hat{g}\hat{H}_{j}^{(nm)}\hat{g}^{-1}$,
then $i=j$ and $\hat{g}\in \hat{H}^{(nm)}_i$.
\end{lem}
%
% Proof
%
\pr
Let $H_i:=H^{(nm)}_i$, $G:=G_{n,m}$ and 
let $c_1,\ldots,c_{nm}$ be free generators of ${\Bbb F}_{nm}$.
For $g\in {\Bbb F}_{nm}$,
let $\hat{g}:=gG\in Q_{n,m}$.
If $i,j$ and $\hat{g}$ satisfy the assumption,
then $[\hat{H}_i, \,
\hat{H}_i\cap  \hat{g}\hat{H}_{j}\hat{g}^{-1}]<\infty$ by definition.
Since $\#\hat{H}_{i}=\infty$,
there exists $\hat{x}\in 
\hat{H}_i\cap  \hat{g}\hat{H}_{j}\hat{g}^{-1}$ such that $\hat{x}\ne 1$.
Then $c_i^lG=xG=gc_j^{l'}g^{-1}G$
for some $l,l'\in {\Bbb Z}$.
For $\hat{\chi}_{u}^{(nm)}$ in  (\ref{eqn:hatchi}),
$u_i^{l}=\hat{\chi}_u^{(nm)}(c_i^lG)
=\hat{\chi}_{u}^{(nm)}(gc_j^{l'}g^{-1}G)
=u_j^{l'}$
for any $u=(u_1,\ldots,u_{nm})\in {\Bbb T}^n\boxtimes {\Bbb T}^m$.
This implies $i=j$ and $l=l'$.
From this, the first statement is verified.
In consequence,
$c_i^lG=gc_i^{l}g^{-1}G$
for some $l\in {\Bbb Z}$.
By the choice of $x$ and $l$, $l\ne 0$.
Since $c_i^lgG=gc_i^{l}G$, $c_i^lg=gc_i^{l}w$ for some $w\in G$.
For $(g',g''):=\phi_{n,m}(g)$,
%
% Equation B.6
%
\begin{equation}
\label{eqn:aib}
(a_{i'}^lg',b_{i''}^lg'')=
\phi_{n,m}(c_i^lg)=
\phi_{n,m}(gc_i^{l} w)=
\phi_{n,m}(gc_i^{l} )=
(g'a_{i'}^l,g''b_{i''}^l)
\end{equation}
where $a_1,\ldots,a_n$
and $b_1,\ldots,b_m$ are free generators of ${\Bbb F}_n$
and ${\Bbb F}_m$, respectively,
and $(i',i'')\in\{1,\ldots,n\}\times\{1,\ldots,m\}$
is defined as $i=m(i'-1)+i''$.
Hence 
$a_{i'}^lg'=g'a_{i'}^l$ and 
$b_{i''}^lg''=g''b_{i''}^l$.
From these and Lemma \ref{lem:multiple}(iii),
$g'\in H_{i'}^{(n)}$ and $g''\in H_{i''}^{(m)}$.
Hence
$\phi_{n,m}(g)=(g',g'')\in H_{i'}^{(n)}\times H_{i''}^{(m)}$.
Therefore $g\in H_{i}G$ by Lemma \ref{lem:hhnm}(ii).
Hence $\hat{g}\in \hat{H}_{i}$.
\qedh

%
% Proposition B.5
% 
\begin{prop}
\label{prop:qnml}
Let  $Q_{n,m}$ and $\hat{H}_{l}^{(nm)}$ 
be as in (\ref{eqn:qnm}) and (\ref{eqn:hath}),
respectively.
\begin{enumerate}
%(i)
\item
For any $l\in\{1,\ldots,nm\}$,
$\lambda_{Q_{n,m}/\hat{H}_l^{(nm)}}$
is irreducible.
%(ii)
\item
For  $l,l'\in\{1,\ldots,nm\}$,
$\lambda_{Q_{n,m}/\hat{H}_l^{(nm)}}
\cong 
\lambda_{Q_{n,m}/\hat{H}_{l'}^{(nm)}}$ if and only if $l=l'$.
\end{enumerate}
\end{prop}
%
% Proof 
%
\pr
(i) 
From Lemma \ref{lem:hij}, ${\rm Com}_{Q_{n,m}}(\hat{H}^{(nm)}_l)
=\hat{H}^{(nm)}_l$.
From this and Theorem \ref{Thm:mackey}(i), the statement holds.

\noindent
(ii)
From Lemma \ref{lem:hij},
$\hat{H}^{(nm)}_l$ and 
$\hat{H}^{(nm)}_{l'}$ are quasiconjugate if and only if 
$l=l'$.
From this and Theorem \ref{Thm:mackey}(ii), the statement holds.
\qedh

\noindent
{\it Proof of Theorem \ref{Thm:ldotss}.}
(i)
From Lemma \ref{lem:equiv},
the statement holds.

\noindent
(ii)
Let $K_l:=G_{n,m}H^{(nm)}_{l}$.
By identifying ${\Bbb F}_{nm}/K_l$ with $Q_{n,m}/\hat{H}_l^{(nm)}$
as a ${\Bbb F}_{nm}$-homogeneous space,
we see that
%
% Equation B.8
%
\begin{equation}
\label{eqn:bbbf}
\lambda_{{\Bbb F}_{nm}/K_l}(g^{(nm)}_l)
=\lambda_{Q_{n,m}/\hat{H}_l^{(nm)}}(\hat{g}^{(nm)}_l)
\quad (l=1,\ldots,nm).
\end{equation}
Hence $\lambda_{{\Bbb F}_{nm}/K_l}({\Bbb F}_{nm})
=\lambda_{Q_{n,m}/\hat{H}_l^{(nm)}}(Q_{n,m})$.
From Proposition \ref{prop:qnml}(i),
$\lambda_{{\Bbb F}_{nm}/K_l}$ is also irreducible.

\noindent
(iii)
From Proposition \ref{prop:qnml}(ii) and
(\ref{eqn:bbbf}),
the statement holds.
\qedh

%%%%%%%%%%%%%%%%%%%%%%%%%%%%%
% 
% Appendix C
%
\sftt{Proof of Proposition \ref{prop:lambdab}(ii)}
\label{section:appthree}
In this section,
we prove Proposition \ref{prop:lambdab}(ii).
For this purpose, we prove the following proposition.
%
% Proposition C.1
% 
\begin{prop}
\label{prop:icc}
Let $Q_{n,m}$ be as in (\ref{eqn:qnm}).
If $n,m\geq 2$,
then the group
$Q_{n,m}$ is ICC, that is,
every conjugacy class in $Q_{n,m}$,
other than its unit is infinite \cite{Blackadar2006}.
\end{prop}
%
% Proof
%
\pr
Rewrite $\phi:=\phi_{n,m}$ and $G:=G_{n,m}$,
and let $\hat{g}:=gG\in Q_{n,m}$ 
and $(g',g''):=\phi(g)$
for $g\in {\Bbb F}_{nm}$.
Fix $\hat{g}\in Q_{n,m}\setminus\{1\}$.
By the choice of $g$,
$(g',g'')\ne (1,1)$.

Assume $g'\ne 1$.
Choose an infinite sequence $\{a_l :l\geq 1\}\subset {\Bbb F}_n$
such that 
$a_lg' a_l^{-1} \ne a_{l'}g' a_{l'}^{-1}$ when $l\ne l'$.
Since ${\Bbb F}_n$ is ICC,
such a sequence always exists.
By the choice of $\{a_l\}$,
$a_l\ne a_{l'}$ when $l\ne l'$.
From Lemma \ref{lem:representation}(i),
there exists $\{(b_l,c_l):l\geq 1 \}\subset {\Bbb F}_{m}\times {\Bbb F}_{nm}$
such that $\phi(c_l)=(a_l,b_l)$  for any $l\geq 1$.
By the choice of $\{c_l\}$,
$\hat{\phi}(c_l gc_{l}^{-1}G)
=(a_lg' a_l^{-1}, b_lg'' b_l^{-1})$
where $\hat{\phi}:Q_{n,m}\to {\Bbb F}_{n}\times {\Bbb F}_{m}$
denotes the natural homomorphism induced by $\phi$.
From this and the choice of $\{a_l\}$,
$\hat{\phi}(c_l gc_{l}^{-1}G)
\ne 
\hat{\phi}(c_{l'} gc_{l'}^{-1}G)$
when $l\ne l'$.
From this,
$\hat{c}_l \hat{g}\hat{c}_{l}^{-1}
=c_l gc_{l}^{-1} G\ne 
c_{l'}gc_{l'}^{-1}G
=
\hat{c}_{l'}\hat{g}\hat{c}_{l'}^{-1}$
when  $l\ne l'$.
Therefore $\{\hat{c}_l \hat{g}\hat{c}_{l}^{-1} :l\geq 1\}$
is an infinite subset of the conjugacy class of $\hat{g}$.

If $g'=1$, then $g''\ne 1$.
In a similar way,  we can construct 
an infinite subset of the conjugacy class of $\hat{g}$ from $g''$
by using Lemma \ref{lem:representation}(ii).
Hence the statement is verified.
\qedh

\noindent
{\it Proof of Proposition \ref{prop:lambdab}(ii).}
By definition,
the group $Q_{n,m}$ in (\ref{eqn:qnm})
acts on $\ell^{2}({\Bbb F}_{nm}/G_{n,m})=\ell^2(Q_{n,m})$
by its left regular representation $\lambda_{Q_{n,m}}$.
For the natural left action $L'$ of $Q_{n,m}$
on $Q_{n,m}$,
$L(g)(hG_{n,m})=ghG_{n,m}=L'(gG_{n,m})(hG_{n,m})$ for any 
$g,h\in {\Bbb F}_{nm}$
where $L$ denotes the natural left action of ${\Bbb F}_{nm}$
on ${\Bbb F}_{nm}/G_{n,m}$.
Hence $\lambda_{{\Bbb F}_{nm}/G_{n,m}}({\Bbb F}_{nm})
=\lambda_{Q_{n,m}}(Q_{n,m})$.
By Proposition \ref{prop:icc},
$\lambda_{Q_{n,m}}$ is a type II$_1$ factor representation
(\cite{Blackadar2006}, III.3.3.7 Proposition). 
Hence the statement holds.
\qedh

%%%%%%%%%%%%%%%%%%%%%%%%%%%%%%%%
%
% Reference 
%

%
\label{Lastpage}

\begin{thebibliography}{99}
%
% AAAAAA
%
\bibitem{Abe}E.\ Abe, 
Hopf algebras.
Cambridge University Press, 1977.
%
\bibitem{AWW}C. Akemann, S. Wassermann,  N. Weaver,
Pure states on free group  C$^*$-algebras.
Glasgow Math.\ J.\ 52 (2010) 151--154.
%
\bibitem{Blackadar2006}B. Blackadar,  
Operator algebras. Theory of C$^{*}$-algebras and
von Neumann algebras. 
Springer-Verlag Berlin Heidelberg New York, 2006.
%
\bibitem{BO}N.P. Brown, N.\ Ozawa,  
C$^{*}$-algebras and finite-dimensional approximations.
American Mathematical Society, 2008.
%
\bibitem{BH}M.\ Burger, P.\ de la Harpe, 
Constructing irreducible representations of 
discrete groups. 
Proc.\ Indian Acad.\ Sci.\ Math.\ Sci.\ 107 (1997)
 223--235.
%
\bibitem{Davidson}K.R. Davidson,
C$^{*}$-algebras by example.
American Mathematical Society, 1996.
%
\bibitem{ES}M. Enock,   J.M. Schwartz, 
Kac algebras and duality of locally compact groups.
Springer-Verlag, 1992.
%
%
%HHHHHHHHH
%
\bibitem{Hatem}H.\ Hatem,
Decomposition of quasi-regular representations
induced from discrete subgroups 
of nilpotent Lie groups.
Lett.\ Math.\ Phys.\ 81 (2007) 135--150.
%
\bibitem{Herschend}M.\ Herschend,  
On the representation ring of a quiver.
Algebr.\ Represent.\ Theor.\ 12 (2009) 513--541.
%
\bibitem{Kassel}C. Kassel, 
Quantum groups. 
Springer-Verlag, 1995.
%
%%%%%%%%%%%%%%%%%%%%%%%
%
\bibitem{TS01}K.  Kawamura, 
A tensor product of representations of Cuntz algebras.
Lett.\ Math.\ Phys.\ 82 (2007) 91--104.
%
\bibitem{TS02}K. Kawamura,  
C$^{*}$-bialgebra defined by the direct sum of Cuntz algebras.
J.\ Algebra 319 (2008) 3935--3959.
%
\bibitem{TS05}K. Kawamura, 
C$^{*}$-bialgebra defined as the direct sum of Cuntz-Krieger algebras.
Comm.\ Algebra 37 (2009) 4065--4078.
%
\bibitem{TS07}K.\ Kawamura, 
Tensor products of type III factor representations of Cuntz-Krieger algebras.
Algebr.\ Represent.\ Theor. 
DOI 10.1007/s10468-012-9362-2,
math.OA/0805.0667v1.
%%%%%%%%%%%%%%%%%%%%%%%%%%%%%%%%%%%%%
%
\bibitem{KV}J. Kustermans,   S.  Vaes,
The operator algebra approach to quantum groups.
Proc. Natl. Acad. Sci. USA 97(2)  (2000) 547--552. 
%
\bibitem{Mackey}G.W. Mackey, 
the theory of unitary group representations.
The University of Chicago Press 
Chicago and London,  1976.
%
\bibitem{MKS}W. Magnus, A. Karrass, D. Solitar,
Combinatorial group theory,
Interscience Publishers, 1966.
%
\bibitem{MNW}T. Masuda,  Y.  Nakagami,   S.L.  Woronowicz, 
A $C^{*}$-algebraic framework for quantum groups.
Int.\ J.\ Math.\ 14 (2003) 903--1001.
%
\bibitem{Ped}G.K. Pedersen, 
$C^{*}$-algebras and their automorphism groups.
Academic Press, 1979.
%
\bibitem{Yoshizawa}H. Yoshizawa,
Some remarks on unitary representations of the free group.
Osaka Math.\ J.\ 3(1) (1951) 55--63.
\end{thebibliography}
\end{document}